\documentclass[3p,sort & compress]{elsarticle}
\usepackage{amssymb}
\usepackage{latexsym}
\usepackage{amsfonts}
\usepackage{amsthm}
\usepackage{amsbsy}
\usepackage{amsmath}
\usepackage{graphics,graphicx}
\usepackage{subfigure}
\usepackage{multirow,multicol}
\usepackage{bm}
\usepackage{bbm}
\usepackage{color}
\usepackage{float}
\usepackage{calc}
\usepackage{caption}
\usepackage{dsfont}
\usepackage{ifthen}
\usepackage{mathrsfs}
\usepackage[colorlinks,linkcolor=blue,citecolor=blue]{hyperref}
\usepackage{epstopdf}
\usepackage{epsfig}
\usepackage{algorithm}
\usepackage{algorithmic}
\usepackage{pifont}
\usepackage{natbib}
\usepackage{overpic}
\usepackage{psfrag}
\usepackage{rotating}
\usepackage{stmaryrd}
\usepackage{verbatim}
\usepackage{setspace}
\usepackage{tikz}

\newdefinition{remark}{Remark}

\newdefinition{method}{Method}
\newdefinition{example}{Example}

\newcommand{\orcid}[1]{\href{https://orcid.org/#1}{\includegraphics[width=8pt]{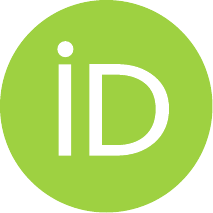}}}

\numberwithin{equation}{section}

\numberwithin{theorem}{section}
\journal{Journal of \LaTeX\ Templates}
\begin{document}
\begin{frontmatter}
    \title{Normalized Fourier-induced PINN method for solving the wave propagation equation in a non-unitized domain over an extended time range}
    \author[label1]{Jichao Ma}\ead{majichao0110@126.com}
    \author[label2]{Dandan Liu}\ead{1486719334@qq.com}
    \author[label3]{Jinran Wu\orcid{0000-0002-2388-3614}}\ead{jinran.wu@uq.edu.au}
    \author[label4,label5]{Xi'an Li\orcid{0000-0002-1509-9328}\corref{cor1}}
    \cortext[cor1]{Corresponding author.}\ead{lixian9131@163.com}
    \address[label1]{School of Mechanical and Electrical Engineering, Qingdao Binhai University, Qingdao 266555, China}
    \address[label2]{College of Arts and Sciences, Qingdao Binhai University, Qingdao 266555, China}
    \address[label3]{School of Mathematics and Physics, The University of Queensland, Brisbane, 4001,Australia}
    \address[label4]{Shandong University, Qingdao 266237, China}
    \address[label5]{Ceyear Technologies Co., Ltd, Qingdao 266555, China}
    \begin{abstract}
    Physics-Informed Neural Networks (PINNs) have gained significant attention for their simplicity and flexibility in engineering and scientific computing. In this study, we introduce a normalized PINN (NPINN) framework to solve a class of wave propagation equations in non-unitized domains over extended time ranges. This is achieved through a normalization technique that involves either spatial or temporal variable normalization. To enhance the capability of NPINN in solving wave equations, we integrate a Fourier-induced deep neural network as the solver, leading to a novel architecture termed NFPINN. Furthermore, we explore different normalization strategies for spatial and temporal variables and identify the optimal normalization approach for our method. To assess the effectiveness and robustness of the proposed NFPINN, we present numerical experiments in both two-dimensional and three-dimensional Euclidean spaces, considering regular and irregular domains. The results confirm the accuracy and stability of our approach.
    \begin{keyword}
    wave propagation; normalization; NFPINN method; large-scale domain; long time range
    % \MSC[2010]65F10 \sep 65F50
    \end{keyword}
    \end{abstract}
\end{frontmatter}

\section{Introduction}\label{sec:intro}
With the continuous advancement of science and technology, the study of the wave propagation problem has acquired crucial significance across various fields, including ocean science and engineering, engineering technology, and military applications, as it provides a fundamental theoretical basis for addressing a wide range of complex problems \cite{liang2019compaction,celep1983spurious}. For instance, in marine engineering, numerous vibration-related issues not only compromise the accuracy and overall performance of power systems but also diminish the stealth capabilities of the system, making it essential to comprehend the fundamental principles governing vibration propagation through the lens of wave theory \cite{mirzajani2018stress}. Mathematically, the general formulation of the wave equation is described in the following:
\begin{equation}\label{eq:wave equation}
    \frac{\partial^2 u(\bm{x}, t)}{\partial t^2} = a^2\Delta u(\bm{x}, t)  + f(\bm{x}, t, u(\bm{x}, t))
\end{equation}
where $\Omega$ denotes a $d-$dimensional polygonal or polyhedral domain in Euclidean space, whose boundary is piecewise Lipschitz and satisfies an interior cone condition. $f(\bm{x}, t, u(\bm{x}, t))$ is a prescribed excitation function and $\Delta $ stands for a standard Laplace operator, i.e.,
\begin{equation*}
    \Delta u(\bm{x}, t)=\sum_{i=1}^{d}\frac{\partial^2 u(\bm{x}, t)}{\partial x_i^2}
\end{equation*}

The solution of Eq.\eqref{eq:wave equation} can be uniquely and high precisely determined through prescribe boundary conditions and initial conditions. Generally, the solution should obey the following boundary constraint
\begin{equation}\label{boundary_condition}
    \mathcal{B}u(\bm{x}, t) = g(\bm{x},t),
\end{equation}
in which, $\mathcal{B}$ is a boundary operator on $\partial \Omega\times(t_0, T]$, such as Dirichlet and Neumman, i.e.,
\begin{equation}\label{Dirichlet_Neumman_bd}
    u(\bm{x}, t) = g_D(\bm{x},t)~~\text{and}~~\frac{\partial u(\bm{x}, t)}{\partial \vec{n}} = g_N(\bm{x},t),
\end{equation}
or both. As well as, the prescribe initial condition is given by
\begin{equation}\label{Initial_condition}
    \mathcal{I}u(\bm{x}, t_0) = h(\bm{x}), ~~\text{for}~~\bm{x}\in\Omega,
\end{equation}
where $\mathcal{I}$ describes the initial operator and $h(\bm{x})$ stands for the initial state of solution, such as 
\begin{equation}\label{Dirichlet_Neumman_initial}
    u(\bm{x}, t_0) = h_K(\bm{x})~~\text{and}~~\frac{\partial u(\bm{x}, t_0)}{\partial t} = h_P(\bm{x}), ~~\text{for}~~\bm{x}\in\Omega.
\end{equation}

The wave propagation equation is a significant partial differential equation widely applied in various fields. It describes different types of waves, including sound waves, light waves, and water waves \cite{uddin2018localized,zhang2018symplectic,sun2017interpolating}. In general, the analytical solution of these partial differential equations poses significant challenges, particularly in complex and inhomogeneous media, thereby necessitating the adoption of numerical methods or asymptotic approaches to obtain approximate solutions \cite{cheng2017fast}. As a result, a wide range of numerical techniques has been developed to address wave-related problems, including but not limited to the finite difference method (FDM) \cite{li2017generalized}, the finite element method (FEM) \cite{li2013fracture}, the boundary element method (BEM), and various meshless methods \cite{liu2006restoring,ma2020meshless,ma2020meshless1,ma2017numerical}, each of which offers distinct advantages depending on the specific problem under consideration.

The finite difference method (FDM) approximates the spatial and temporal derivatives of the wave field function through corresponding spatial and temporal differences \cite{takekawa2018mesh}, thereby enabling high computational efficiency, low memory consumption, and superior accuracy; however, its dependence on a single computational mesh inherently leads to poor numerical stability. In contrast, the finite element method (FEM), which is founded on element-based analysis and relies on a predefined mesh interconnected by nodes, proves particularly suitable for handling complex geometries and addressing a diverse array of scientific and engineering problems \cite{gao2019combining,chen2013combining}. Nevertheless, when dealing with intricate structures, the mesh generation process can be time-consuming, and as the number of mesh nodes increases, computational costs grow significantly.  
(BEM), also known as the boundary integral equation method (BIE), represents an advanced numerical technique that was developed subsequent to the finite element method (FEM), offering distinct advantages in solving boundary-value problems with reduced dimensionality \cite{mantegh2010path,xie2014direct}. It possesses several advantages, including dimensionality reduction and the ability to handle infinite domains, as it only requires discretization of the computational boundary. However, BEM also has limitations, such as the asymmetry of the coefficient matrix in its equations and the necessity of fundamental solutions, which can be difficult to derive for certain complex problems.  
The meshless method, also called the element-free method, represents the interior and boundary of a problem domain using a set of distributed nodes \cite{cheng2014novel,cheng2016analyzing,liu2008numerical}. These nodes, referred to as field nodes, eliminate the necessity for predefined node connectivity when formulating interpolation or approximation expressions of unknown functions, thereby simplifying computational implementation and enhancing flexibility in numerical analysis \cite{liu2017density,vzilinskas2010similarities,lin2014almost,liu2005modeling}. In comparison to conventional finite element (FEM) and finite difference methods (FDM), the meshless method effectively overcomes challenges associated with mesh generation, distortion, and movement, thus offering an innovative and highly efficient approach to engineering analysis \cite{shi2013extended,liu2018improved,gao2017complex,tatari2011finite,darani2017direct}.However, its computational results can be unstable under varying parameters, such as shape parameters, penalty factors, and kernel functions.

Deep neural networks (DNN), as a powerful meshfree approach that eliminates the need for domain discretization, have increasingly attracted the attention of researchers for numerically solving ordinary and partial differential equations (PDEs) as well as inverse problems involving complex geometries and high-dimensional scenarios \cite{weinan2018deep,sirignano2018dgm,chen2018neural,raissi2019physics,khoo2019switchnet,zang2020weak,lyu2022mim}, a trend largely driven by their remarkable universal approximation capability \cite{hauptmann2020deep}. Among these methods, physics-informed neural networks (PINN), whose origins can be traced back to the early 1990s, have once again garnered widespread interest and achieved significant breakthroughs in approximating PDE solutions by embedding physical laws into neural networks, benefiting from the rapid advancements in computer science and technology \cite{raissi2019physics,dissanayake1994neural}. This innovative approach constructs a loss function by integrating the residuals of governing equations with discrepancies in boundary and initial conditions, enabling efficient optimization through automatic differentiation within DNN. To further enhance the performance of PINN, extensive research has focused primarily on two key aspects: refining the selection of residual terms and developing strategies to effectively enforce initial and boundary constraints. In terms of residual terms, various improved formulations such as XPINN \cite{jagtap2020extended}, cPINN \cite{jagtap2020conservative}, two-stage PINN \cite{lin2022two}, and gPINN \cite{yu2022gradient} have been proposed. Furthermore, by cleverly imposing initial and boundary constraints directly within DNN, PINN not only become easier to train but also exhibit reduced computational complexity and achieve highly accurate PDE solutions under complex boundary conditions \cite{berg2018unified,sun2020surrogate,lu2021physics}. Inspired by the order reduction techniques traditionally employed in numerical methods, researchers have explored reformulating high-order PDEs as first-order systems to mitigate the computational difficulties associated with high-order derivatives in DNN, leading to the development of notable approaches such as the deep mixed residual method \cite{lyu2022mim}, the local deep learning method \cite{Yang2021local}, and the deep FOSLS method \cite{cai2020deep,bersetche2023deep}.

To approximate the solution of the wave equation and full waveform inversions (FWIs), Majid Rasht-Behesht et al. applied the PINN to solve the acoustic wave equation and test the presented method with both forward models and FWI case studies \cite{rasht2022physics}. 
For high-dimensional wave equations, Shaikhah Alkhadhr et al. utilized the PINN method to simulate a linear wave equation with a single time-dependent sinusoidal source function e.g.: $\sin (\pi t)$ which is one of the most fundamental modeling equations in medical ultrasound applications, then validating its efficacy and testing its performance \cite{alkhadhr2021modeling}.  
In order to make the diffusion wave model (DWM) widely applied in water conservancy, hydrology and irrigation engineering. Qingzhi Hou et al. applied PINN with novel improvements to solve the DWM for both forward and inverse problems, and obtained remarkable results \cite{hou2024physics}.
To tackle wave propagation in fully nonlinear potential flows with the free surface, Haocheng Lu et al. proposed fully nonlinear free surface physics-informed neural networks (FNFS-PINNs) which offers an approach to address the complexities of modeling nonlinear free surface flows and broaden the scope of PINN for various wave propagation scenarios \cite{lu2024physics}. 
Katayoun Eshkofti et al. utilized PINN to study thermoelastic wave propagation and the Moore-Gibson-Thompson (MGT) coupled thermoelasticity in porous media. They formulated the governing equations for coupled thermoelasticity in a porous half-space based on the MGT heat conduction model, incorporating the thermal relaxation coefficient and strain relaxation factor. Mechanical and thermal shock loading boundary conditions were applied. The PINN method was employed to analyze the behavior of a porous magnesium body, yielding highly accurate results for the system of coupled PDEs\cite{eshkofti2024modified}. Because analytical acoustic wave equation (AWE) solutions rarely exist for complex heterogeneous media, Harpreet Sethi et al. used the mesh-free PINN method to generate AWE solutions, in which  a Fourier neural network (FNN) was introduced to alleviate the spectral bias commonly observed when using fully connected neural network in the conventional PINN approach \cite{sethi2023hard}. 

While many researchers have studied the capacity of PINN for solving wave equations in unit and short time range, studies on its solution in large-scale domains and long time range are rare. In this paper, we develop a normalized Fourier induced PINN method (dubbed NFPINN) to solve wave equations in large-scale domains and long time range. By applying the normalization technique for spatial or temporal variable, then the solver of PINN will easily to handle these space-time inputs and product a satisfactory results for this wave paopagation problem. Moreover, by incorporating a Fourier feature mapping, the ability to efficiently capture the target function and conveniently express its derivatives is significantly enhanced, thereby substantially improving the capacity of the PINN model in solving wave equations.

In summary, the principal contributions of this paper can be outlined as follows:

1. We introduce a novel neural network approach that integrates the PINN method with Fourier feature mapping technology within a subnetwork structure, thereby enhancing the capability to effectively address the wave equation.

2. Through numerical experiments, we demonstrate that the classical PINN method strengthened by a Fourier-induced DNN  still remains inadequate in delivering accurate solutions for the wave equation, particularly in a non-unitized domain and over extended time ranges.

3. We demonstrate the outstanding performance of FPINN in solving a class of wave propagation equation with Dirichlet and Neumann boundary conditions and initial conditions across various dimensional spaces, within a non-unitized domain, and over an extended time range.

The remainder of this paper is structured as follows. In Section \ref{sec:2}, we provide a concise introduction to the fundamental concepts and formulations of PINN, along with its Fourier-enhanced variant. Section \ref{sec:3} presents a unified framework for the NFPINN, designed to solve the wave equation in non-unitized domains and over extended time ranges, while also offering an error analysis of the proposed approach. In Section \ref{sec:4}, we elaborate on the NFPINN algorithm for approximating the solution to the wave equation. Section \ref{sec:5} showcases several wave equation scenarios to assess the feasibility and effectiveness of our proposed method. Finally, Section \ref{sec:6} concludes the paper with a summary of key findings.

\section{Fourier induced PINN to the wave equation}\label{sec:2}
\subsection{Fourier induced deep neural network}\label{sec:FPINN}
In order to approximate the solution of the wave equation governed by Eq.\eqref{eq:wave equation} on interested domain, we would like to design PINN method to infer the all system states of interest under prescribed boundary and initial conditions. The basic solver of PINN framework is generally configured as a DNN model, and the DNN defines the following mapping
\begin{equation}
   \mathcal{F}: \bm{x}\in\mathbb{R}^{d}\Longrightarrow \bm{y}=\mathcal{F}(\bm{x})\in\mathbb{R}^{c}
\end{equation}
with $d$ and $c$ denoting the input and output dimensions, respectively. The DNN function $\mathcal{F}$ is structured as a hierarchical composition of consecutive linear transformations and nonlinear activation functions, following a sequential arrangement.
\begin{equation}
\begin{cases}
\bm{y}^{[0]} = \bm{x}\\
\bm{y}^{[\ell]} = \sigma\circ(\bm{W}^{[\ell]}\bm{y}^{[\ell-1]}+\bm{b}^{[\ell]}), ~~\text{for}~~\ell =1, 2, 3, \cdots\cdots, L
\end{cases}
\end{equation}
where $\bm{W}^{[\ell]} \in  \mathbb{R}^{n_{\ell+1}\times n_{\ell}}, \bm{b}^{[\ell]}\in\mathbb{R}^{n_{\ell+1}}$ are the weights and biases of $\ell$-th hidden layer, respectively, where $n_0=d$ and $n_{L+1}$ corresponds to the output dimension. The symbol $``\circ"$ denotes the element-wise operation, while the function $\sigma(\cdot)$ represents an element-wise activation function. We define the output of a DNN as $\bm{y}(\bm{x};\bm{\theta})$, where $\bm{\theta}$ denotes the set of parameters, including $\bm{W}^{[1]},\cdots \bm{W}^{[L]}$, and, $\bm{b}^{[1]},\cdots \bm{b}^{[L]}$.

\begin{figure}[H]
	\centering
	\includegraphics[scale=0.6]{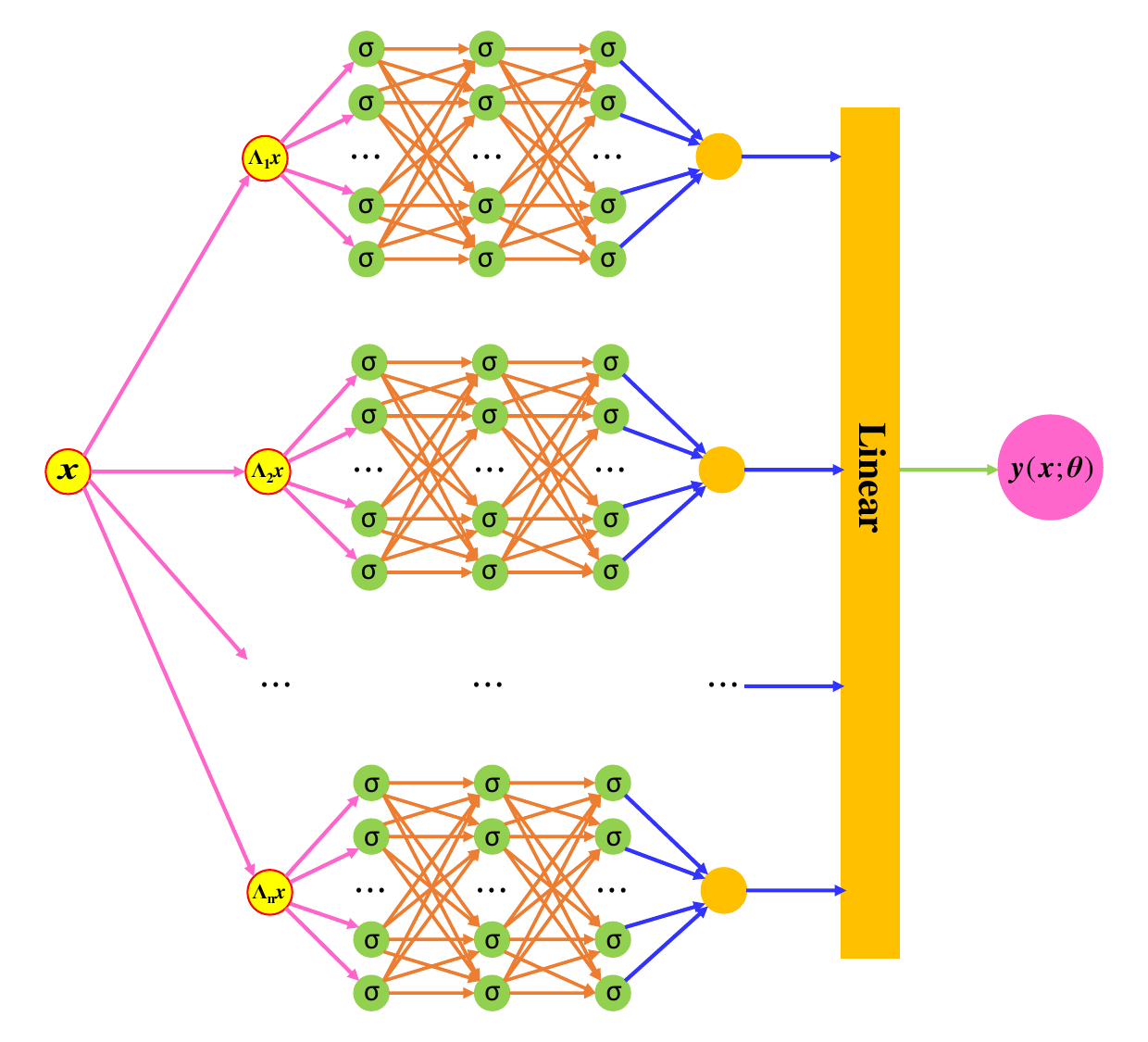}
	\caption{Schematic diagram of FFM-based DNN with $n$ subnetworks, $\sigma$ stands for the activation function}
	\label{Fig:FPINN network}
\end{figure}

To enhance the capacity of a general PINN, a Fourier feature mapping (FFM), composed of sine and cosine functions, is introduced as the activation functions for the first hidden layer in a DNN, as this approach not only alleviates the issue of spectral bias but also enables the network to learn high-frequency components more effectively \cite{rahaman2018spectral, wang2020eigenvector, Matthew2020Fourier, Xu_2020, LI2023114963}, thereby improving its overall performance. For convenience, the PINN model with a classical DNN architecture induced by FFM is called as Fourier PINN (abbreviated as FPINN). As we have known, the FFM composed by sine and cosine functions is a global basis function and is optimal, leading to the following formulation.
\begin{equation}\label{fourier}
\zeta(\bm{x}) = 
\left[\begin{array}{c}
\cos(\bm{\Lambda} \bm{x})\\
\sin(\bm{\Lambda} \bm{x}) 
\end{array}
\right]
\end{equation}
where $\bm{\Lambda}$ serves as a pre-input vector or matrix that can either be trainable or non-trainable, plays a crucial role in directly influencing the first hidden layer of the DNN by ensuring consistency with the number of nodes in that layer, as illustrated in the schematic diagram of the FFM-based DNN with $n$ subnetworks shown in Fig.\ref{Fig:FPINN network}.

\subsection{Optimizing the loss function to solve the wave equation}
Suppose the surrogate model, denoted by $u_{NN}(\bm{x}, t)$, expressed by DNN model predicted, then optimal approximated solution corresponding to the wave equation can be derived through the process of minimizing the following loss function
\begin{equation}\label{loss2PINN}
Loss=Loss_{PDE} + Loss_{BC} + Loss_{IC}
\end{equation}
with
\begin{equation}\label{subloss2PINN}
\begin{aligned}
&Loss_{PDE} = \frac{1}{N_R}\sum_{i=1}^{N_R}\left| \frac{\partial^2 u_{NN}(\bm{x}_R^{i}, t_R^{i})}{\partial t^2} - a^2\Delta u_{NN}(\bm{x}_R^{i}, t_R^{i}) - f(\bm{x}_R^{i}, t_R^{i}, u_{NN}(\bm{x}_R^{i}, t_R^{i}))\right|^2 \\
&Loss_{BC} = \frac{1}{N_B}\sum_{i=1}^{N_B}\bigg{|}\mathcal{B}u_{NN}\left(\bm{x}_B^i, t_B^i\right)-g(\bm{x}_B^i, t_B^i)\bigg{|}^2\\
&Loss_{IC} = \frac{1}{N_I}\sum_{i=1}^{N_I}\bigg{|}\mathcal{I}u_{NN}\left(\bm{x}_I^i,t_0\right)-h(\bm{x}_I)\bigg{|}^2
\end{aligned}
\end{equation}
where $\gamma$ is a weighting parameter for the boundary loss. Here, $Loss_{PDE}$, $Loss_{BC}$ and $Loss_{IC}$ denote the residual errors for the governing differential equations (PDE),  the specified boundary conditions (BC) and initial conditions (IC), respectively. In the case where additional data are accessible within the domain, it becomes possible to incorporate an extra loss term, which serves to quantify and account for the discrepancy between the approximated solution and the raw data, thereby enhancing the overall accuracy of the model.
\begin{equation}
    Loss_{D} = \frac{1}{N_D}\sum_{i=1}^{N_D}\bigg{|}u_{NN}(\bm{x}_D^i,t_D^i)-u_{Data}^i\bigg{|}^2
\end{equation}

\begin{figure}[H]
    \centering
    \includegraphics[scale=0.4]{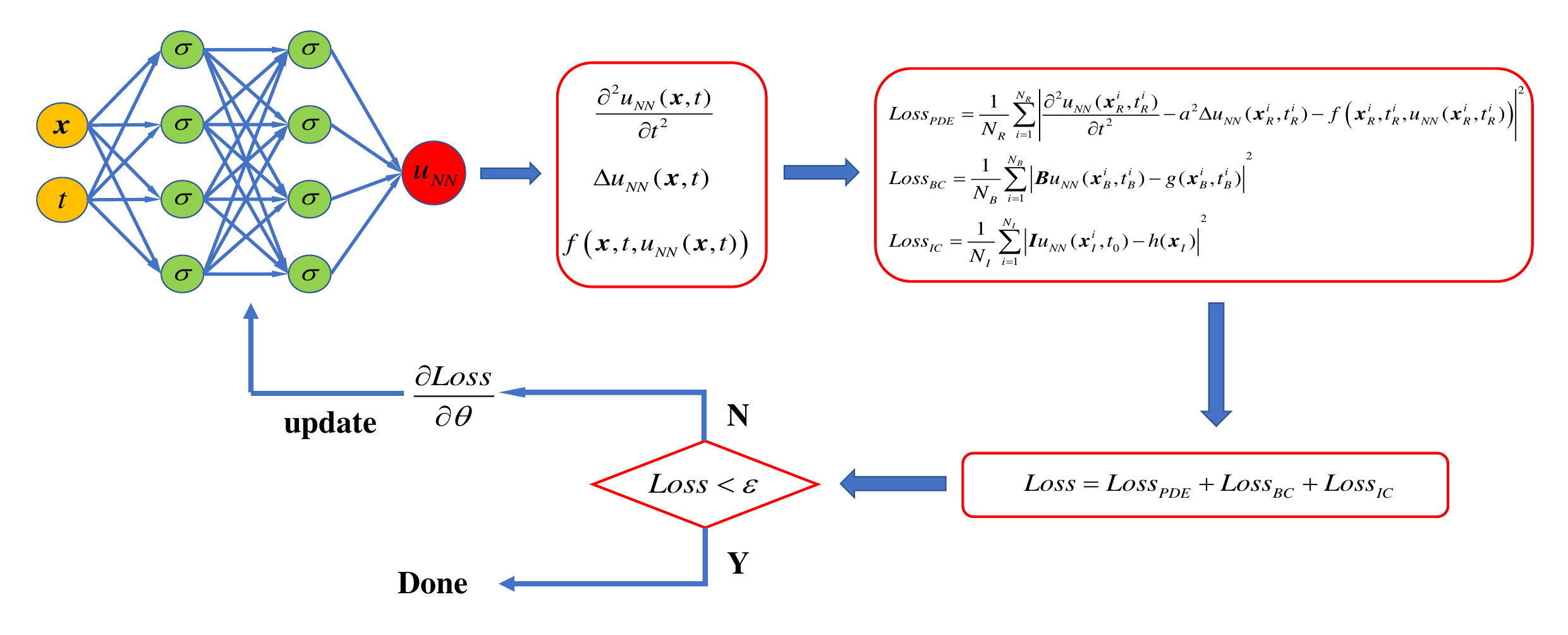}
    \caption{FPINN framework to solve wave propagation}
    \label{Fig:FPINN framework}
\end{figure}

\section{Normalized Fourier induced PINN to solve the wave equation}\label{sec:3}
\subsection{The analysis of general PINN and FPINN method to the wave equation in two different scale range}

Although various PINN models have been successfully applied to the study of ordinary and partial differential equations, particularly in the case of the wave equation, our investigation shows that their performance deteriorates in large scale domain and long time range, potentially leading to non-convergence.

For example, let us consider two scenarios for two-dimensional wave propagation equation with Dirichlet boundary in $\Omega_1=[0,2\pi]\times[0,2\pi], t\in(0,2)$. $\Omega_2=[0,10\pi]\times[0,10\pi], t\in(0,10)$, respectively. The governed equation is  

\begin{equation}\label{eq:Example_1_govern}
    \frac{\partial^2 u}{\partial t^2} = \frac{1}{2}\bigg{(}\frac{\partial^2 u}{\partial x_1^2} + \frac{\partial^2 u}{\partial x_2^2}\bigg{)} + 12t^2
\end{equation}

An analytical solution is given by
\begin{equation*}
    u(x_1,x_2,t) = t^4 + \sin(x_1)\cdot\sin(x_2)\cdot\sin(t).
\end{equation*}

Since the boundary and initial constraint functions can be directly derived from the exact solution, we will not explicitly state them here.

In this experiment, the solvers for PINN and FPINN are configured as a DNN and a FFM-based DNN with $N$ subnetworks, respectively, and the scale factors are set as (1, 2, 3, 4, 5, 6, 7, 8, 9, 10), and each subnetwork is configured with sizes of (20, 15, 15, 10). The first hidden layer of all subnetworks employs Fourier feature mapping as the activation function (see Eq.\eqref{fourier}), while the activation functions for the other layers (except for the output layer) are selected as $\text{GELU}(x) = x \cdot \frac{1}{2}[1+erf(\frac{x}{\sqrt{2}})]$, where $erf(x)$ is Gaussian error function, and the output layers of all subnetworks are linear. We train the previously mentioned PINN and FPINN models for 30,000 epochs, performing testing every 1,000 epochs during the training process. The optimizer is set as Adam with an initial learning rate of 0.01 and the learning rate will decay by 3.5\% for every 100 epochs.

\begin{figure}[H]
    \centering
    \subfigure[Analytical solution]{
        \label{E1-01} 
        \includegraphics[scale=0.45]{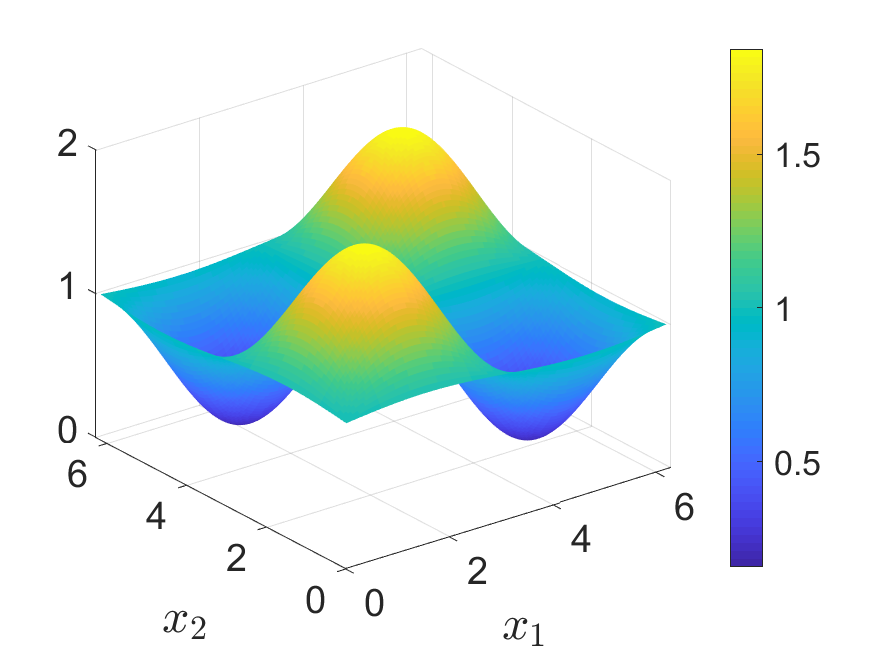}
    }
    \subfigure[Predicted solution for PINN]{
        \label{E1-02} 
        \includegraphics[scale=0.45]{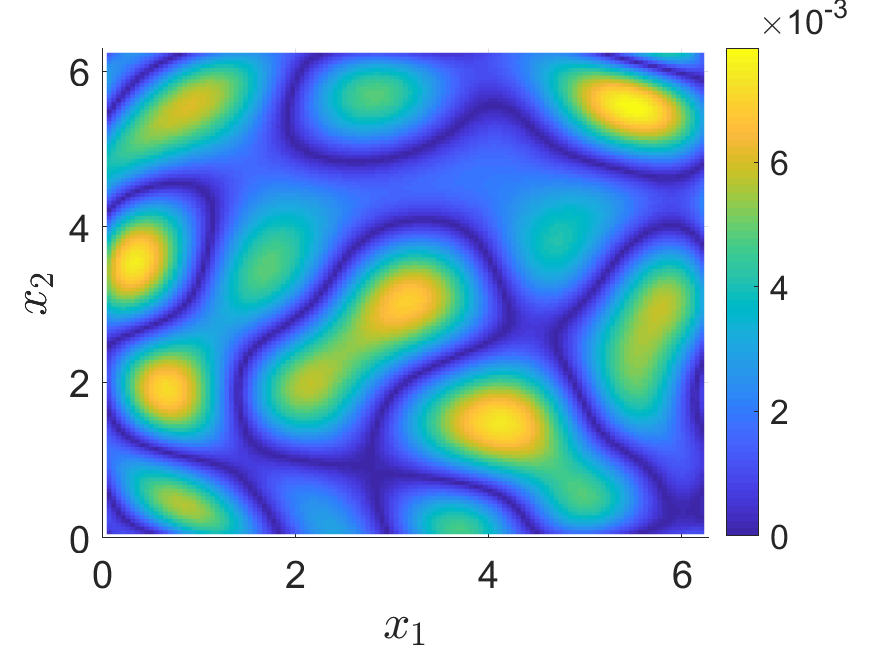}
    }
    \subfigure[Predicted solution for FPINN]{
        \label{E1-03} 
        \includegraphics[scale=0.45]{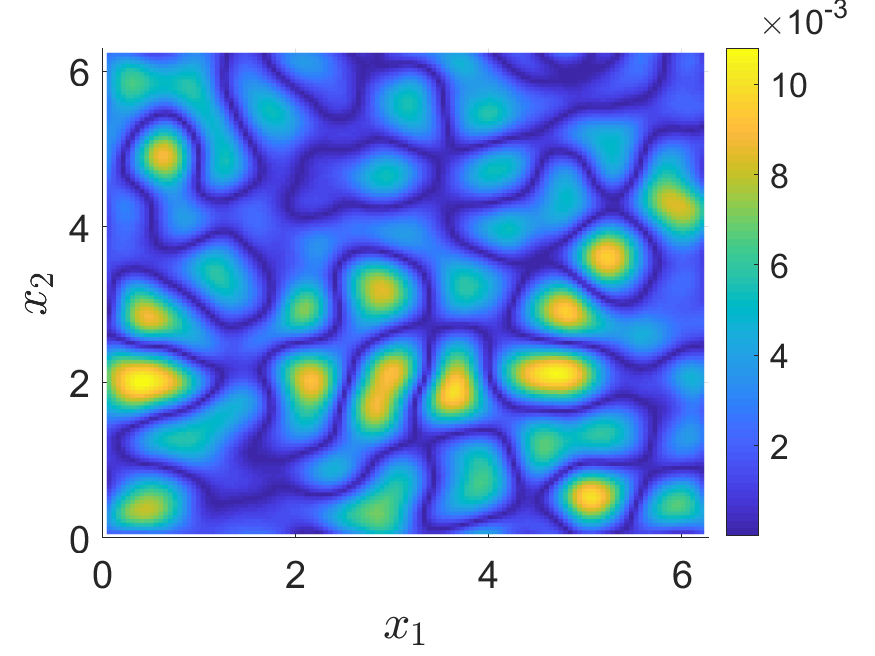}
    }
    \subfigure[PINN vs FPINN ]{
        \label{E1-04}
        \includegraphics[scale=0.425]{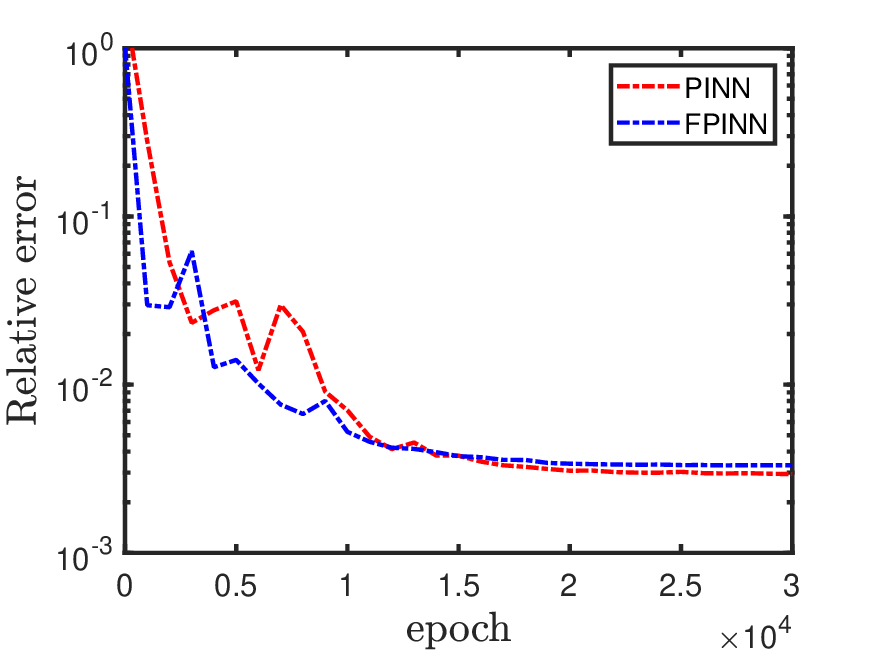}
    }
    \caption{Numerical results of PINN and FPINN methods for test example in domain $\Omega_1$.}
    \label{Fig:Example_1_small}
\end{figure}

\begin{figure}%[H]
    \centering
    \subfigure[Analytical solution]{
        \label{E1-05} 
        \includegraphics[scale=0.45]{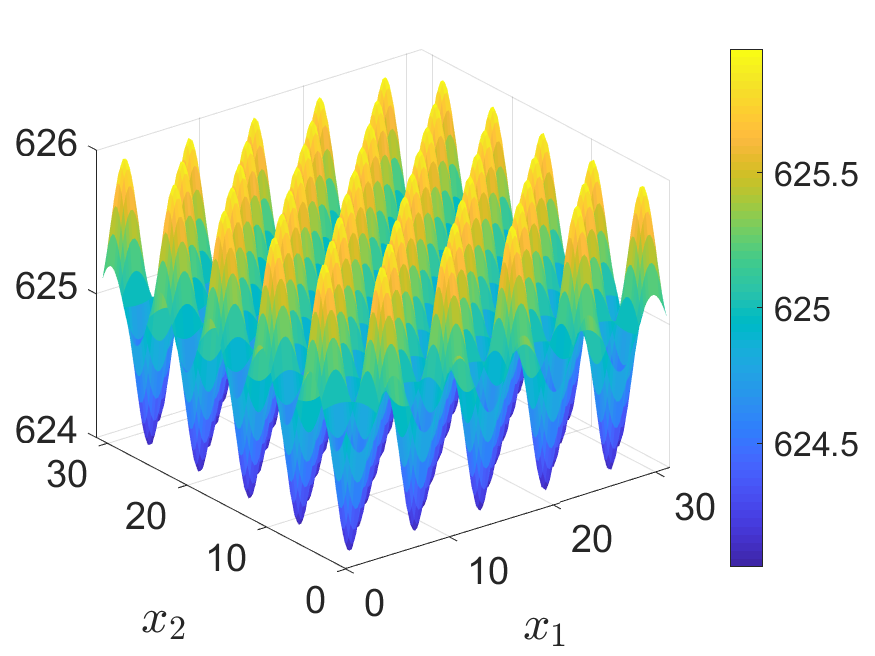}
    }
    \subfigure[Predicted solution for PINN]{
        \label{E1-06} 
        \includegraphics[scale=0.45]{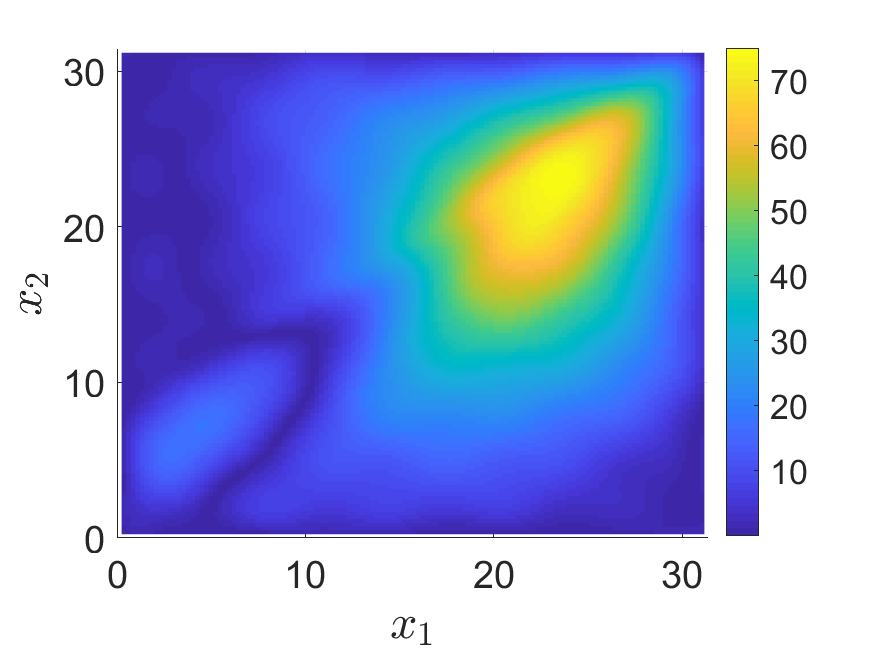}
    }
    \subfigure[Predicted solution for FPINN]{
        \label{E1-07} 
        \includegraphics[scale=0.45]{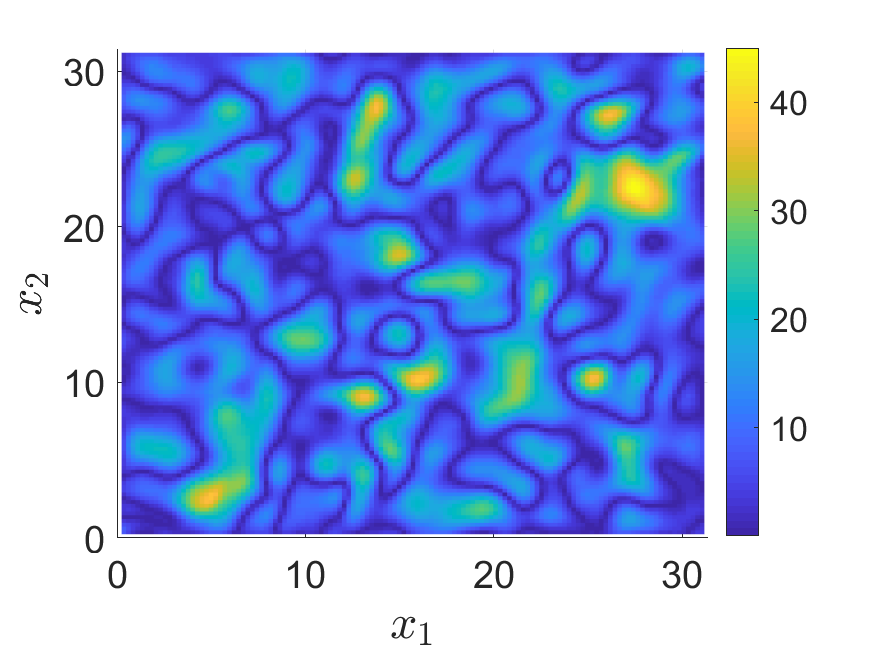}
    }
    \subfigure[PINN vs FPINN ]{
        \label{E1-08}
        \includegraphics[scale=0.425]{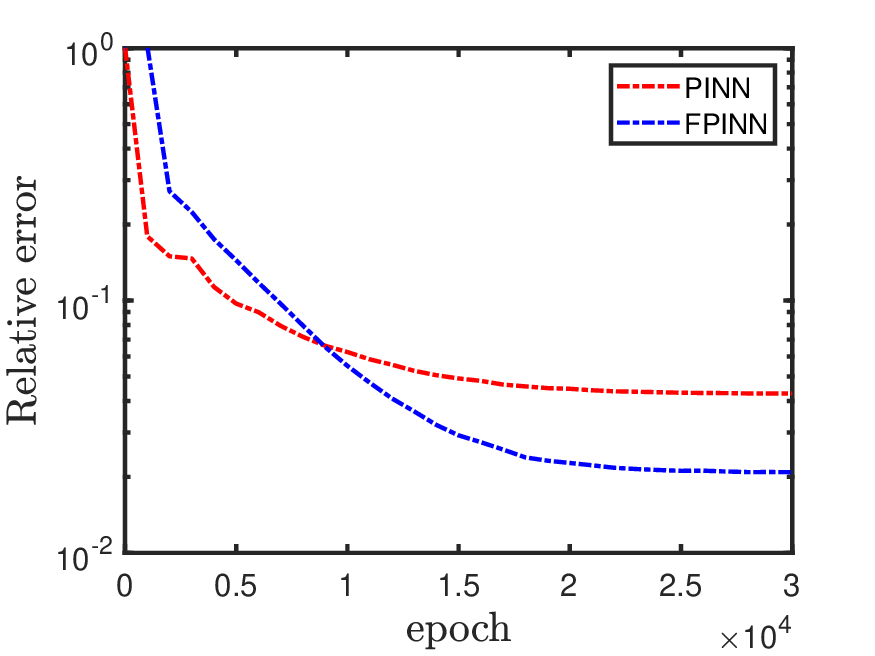}
    }
    \caption{Numerical results of PINN and FPINN methods for test example in domain $\Omega_2$.}
    \label{Fig:Example_1_large}
\end{figure}

By employing Latin Hypercube Sampling (LHS) to strategically select training points across the temporal and spatial domains, we ensure that, in each epoch, the model is trained with $N_R$=1500 collocation points for the governed equation, $N_B$ =300 boundary points for the boundary condition, and $N_I$=700 initial points for the initial condition, thereby facilitating comprehensive learning that enables the model to effectively capture intricate spatiotemporal patterns. In addition, 16384 collocations sampled from regular domain $\Omega_1$ at $t = 0.5$ and $\Omega_2$ at $t = 2.5$ are used as the testing set to validate the feasibility of model. The results are plotted in Fig.\ref{Fig:Example_1_small} and Fig.\ref{Fig:Example_1_large}.

From Fig.\ref{Fig:Example_1_small} and Fig.\ref{Fig:Example_1_large}, the predicted solutions by PINN and FPINN methods are pretty close to analytical solution of Eq.\eqref{eq:Example_1_govern} in small-scale domain and short-time range. However, PINN and FPINN methods cannot approximate well the solution of Eq.\eqref{eq:Example_1_govern} in large-scale domain and long-term range,even though the approximation producted by FPINN method is better than that of PINN method. Moreover, we apply the aforementioned PINN and FPINN methods with various hyperparameter configurations, including different learning rates, numbers of subnetworks, and network sizes, yet we still fail to achieve a satisfactory result.

\subsection{NFPINN method to solve wave propagation in large-scale domain and long-time range}\label{sec:FMPINN}
As described in the above, the FPINN method is elaborated to deal with the wave equation Eq.~\eqref{eq:wave equation} with Dirichlet or Neumman boundary condition. Generally, this problem is solved in a small-scale domain and short-time range, but the solutions for many application scenarios are required in a large-scale domain and long-time range. In this section, we introduce the novel concept of the NFPINN method, which incorporates skilly spatial and/or  temporal normalization strategies with the FPINN method.

\subsubsection{Spatial NFPINN (S-NFPINN)}
First, we normalize the spatial domain $\Omega^d$ into $\overline{\Omega} = [0,1]^d$ for each component of independent variable
\begin{equation*}
    \tilde{x}_i = \frac{x_i-x_i^{min}}{s_i}
\end{equation*}
where the scale factor $s_i = x_i^{max} - x_i^{min}$, where $x_i^{max}$ and $x_i^{min}$ stand for the maximum and minimum values for component $x_i$ of $\bm{x}$, respectively. The temporal domain is expressed as $T=(t_0, t_{max}]$. Then, Eq.\eqref{eq:wave equation} can be rewritten as 

\begin{equation}\label{normalized:wave_space}
    \frac{\partial^2 \tilde{u}(\tilde{\bm{x}}, t)}{\partial t^2} = a^2\bigg{(}\sum_{i=1}^{d}\frac{1}{s^2_i}\frac{\partial^2 \tilde{u}(\tilde{\bm{x}}, t)}{\partial \tilde{x}_i^2} \bigg{)} + \tilde{f}\big{(}\tilde{\bm{x}}, t, \tilde{u}(\tilde{\bm{x}}, t)\big{)}, ~~\text{for}~~ (\tilde{\bm{x}}, t)\in(\overline{\Omega}, T)
\end{equation}
with 
\begin{equation*}
    \tilde{f}(\tilde{\bm{x}}, t, \tilde{u}(\tilde{\bm{x}}, t)) = f(x_1^{min}+s_1\tilde{x}_1, x_2^{min}+s_2\tilde{x}_2, \cdots, x_d^{min}+s_d\tilde{x}_d, t, \tilde{u}(\tilde{x}_1, \tilde{x}_2, \cdots, \tilde{x}_d, t)).
\end{equation*}

The corresponding Dirichlet type and Neumman type boundary condition from Eq.\eqref{Dirichlet_Neumman_bd} is written as
\begin{equation}\label{normalized:wave_space:Dirichlet}
    \tilde{u}(\tilde{\bm{x}},t) = \tilde{g}_D(\tilde{\bm{x}},t)
\end{equation}
and
\begin{equation}\label{normalized:wave_space:Neumman}
    \frac{\partial \tilde{u}(\tilde{\bm{x}},t)}{\partial \vec{n}} = \frac{1}{s_1}\frac{\partial \tilde{u}(\tilde{\bm{x}},t)}{\partial \tilde{x}_1}\cdot n_1 + \frac{1}{s_2}\frac{\partial \tilde{u}(\tilde{\bm{x}},t)}{\partial \tilde{x}_2}\cdot n_2 + \cdots + \frac{1}{s_d}\frac{\partial \tilde{u}(\tilde{\bm{x}},t)}{\partial \tilde{x}_d}\cdot n_d = \tilde{g}_N(\tilde{\bm{x}},t),
\end{equation}
respectively. Similarly, the initial condition Eq.\eqref{Dirichlet_Neumman_initial} can be modified as

\begin{equation}\label{normalized:wave_space:initial}
    \tilde{u}(\tilde{\bm{x}},t_0) = \tilde{h}_K(\tilde{\bm{x}}), ~~\text{and}~~ \frac{\partial \tilde{u}(\tilde{\bm{x}},t_0)}{\partial t} = \tilde{h}_P(\tilde{\bm{x}}).
\end{equation}

\subsubsection{Temporal NFPINN (T-NFPINN)}
Second, casting the original time range into the interval $\overline{T}=[0,1]$ by time normalization technique, we adeptly mitigate challenges attributed to protracted temporal dependencies, i.e.,
\begin{equation*}
    \tilde{t} = \frac{t-t_0}{s_T}
\end{equation*}
with $s_T = t_{max}-t_0$. Then, we have the following normalization formulation
\begin{equation}\label{normalized:wave_time}
    \frac{1}{s_T^2}\frac{\partial^2 \tilde{u}(\bm{x},\tilde{t})}{\partial \tilde{t}^2} = a^2\bigg{(}\sum_{i=1}^{d}\frac{\partial^2 \tilde{u}(\bm{x},\tilde{t})}{\partial x_i^2} \bigg{)} + \tilde{f}\big{(}\bm{x},\tilde{t},\tilde{u}(x,\tilde{t})\big{)}, ~~\text{for}~~ ({\bm{x}},\tilde{t})\in({\Omega},\overline{T})
\end{equation}
with 
\begin{equation*}
    \tilde{f}(\bm{x}, \tilde{t}, \tilde{u}(\bm{x}, \tilde{t})) = f(\bm{x}, t_0+s_T\tilde{t},\tilde{u}(\bm{x}, \tilde{t})).
\end{equation*}

The corresponding Dirichlet type and Neumman type boundary condition from Eq.\eqref{Dirichlet_Neumman_bd} is written as
\begin{equation}
    \tilde{u}(\bm{x},\tilde{t}) = \tilde{g}_D(\bm{x},\tilde{t})
\end{equation}
and
\begin{equation}
    \frac{\partial \tilde{u}(\bm{x}, \tilde{t})}{\partial \vec{n}} = \tilde{g}_N(\bm{x}, \tilde{t}),
\end{equation}
respectively. Similarly, the initial condition Eq.\eqref{Dirichlet_Neumman_initial} can be modified as

\begin{equation}
    \tilde{u}(\bm{x}, t_0) = \tilde{h}_K(\bm{x}), ~~\text{and}~~ \frac{1}{s_T}\frac{\partial \tilde{u}(\bm{x}, t_0)}{\partial \tilde{t}} = \tilde{h}_P(\bm{x}).
\end{equation}

\subsubsection{Spatial and temporal NFPINN (ST-NFPINN)}
As the final strategy, the spatial and temporal normalized formulation is applied to Eq.\eqref{eq:wave equation}, thereby enabling the governing equation of Eq.\eqref{eq:wave equation} to be reformulated in terms of the scaled and characteristic quantities.

\begin{equation}\label{normalized:wave_time_space}
    \frac{1}{s_T}\frac{\partial^2 \tilde{u}(\tilde{\bm{x}},\tilde{t})}{\partial \tilde{t}^2} = a^2\bigg{(}\sum_{i=1}^{d}\frac{1}{s^2_i}\frac{\partial^2 \tilde{u}(\tilde{\bm{x}},\tilde{t})}{\partial \tilde{x}_i^2} \bigg{)} + \tilde{f}\big{(}\tilde{\bm{x}},\tilde{t},\tilde{u}(\tilde{\bm{x}},\tilde{t})\big{)}, ~~\text{for}~~ (\tilde{\bm{x}},\tilde{t})\in(\overline{\Omega},\overline{T})
\end{equation}
in which, the definition for $\overline{\Omega}$, $\overline{T}$, $s_i$ and $s_T$ as same as the above descriptions, and $\tilde{f}\big{(}\tilde{\bm{x}},\tilde{t},\tilde{u}(\tilde{\bm{x}},\tilde{t})\big{)}$ is easy to obtain. The corresponding boundary conditions Eq.\eqref{boundary_condition} on unit space domain and time range can be expressed as $\mathcal{B}\tilde{u}(\tilde{\bm{x}}, \tilde{t}) = \tilde{g}(\tilde{\bm{x}},\tilde{t})$, in particular, the Neumman type for boundary constraint is written as
\begin{equation}
    \tilde{u}(\tilde{\bm{x}},\tilde{t}) = \tilde{g}_D(\tilde{\bm{x}},\tilde{t})
\end{equation}
and
\begin{equation}
    \frac{\partial \tilde{u}(\tilde{\bm{x}},\tilde{t})}{\partial \vec{n}} = \frac{1}{s_1}\frac{\partial \tilde{u}(\tilde{\bm{x}},\tilde{t})}{\partial \tilde{x}_1}\cdot n_1 + \frac{1}{s_2}\frac{\partial \tilde{u}(\tilde{\bm{x}},\tilde{t})}{\partial \tilde{x}_2}\cdot n_2 +\cdots+ \frac{1}{s_d}\frac{\partial \tilde{u}(\tilde{\bm{x}},\tilde{t})}{\partial \tilde{x}_d}\cdot n_d = \tilde{g}_N(\tilde{\bm{x}},\tilde{t}),
\end{equation}
respectively. Similarly, the initial condition Eq.\eqref{Dirichlet_Neumman_initial} can be modified as
\begin{equation}
    \tilde{u}(\tilde{\bm{x}}, t_0) = \tilde{h}_K(\tilde{\bm{x}}), ~~\text{and}~~ \frac{1}{s_T}\frac{\partial \tilde{u}(\tilde{\bm{x}}, t_0)}{\partial \tilde{t}} = \tilde{h}_P (\tilde{\bm{x}}).
\end{equation}

\subsection{The loss function of NFPINN and its optimization}
Building upon the aforementioned discussions, the procedure of the NFPINN algorithm, which is designed to address the wave equation Eq.\eqref{eq:wave equation} within finite-dimensional spaces, is outlined in the following.

During the $k_{th}$ iteration step, a set of sampled collocation points $\mathcal{S}^k=\mathcal{S}_R^k\cup \mathcal{S}_B^k\cup\mathcal{S}_I^k$ is given, where $\mathcal{S}_R^k$, $\mathcal{S}_B^k$, $\mathcal{S}_I^k$ are the points set sampled from interior, boundary and initial, respectively. Then the loss function of S-NFPINN method for the wave problem with Dirichlet boundary can be obtained by Eqs. \eqref{normalized:wave_space}, \eqref{normalized:wave_space:Dirichlet} and \eqref{normalized:wave_space:initial}, it is expressed as follows
\begin{equation}\label{Loss_NFPINN_Dirichlet}
\begin{aligned}
    Loss_{Dir}&=\frac{1}{N_R}\sum_{i=1}^{N_R}\left| \frac{\partial^2 \tilde{u}_{NN}(\tilde{\bm{x}}_R^i,t_R^i)}{\partial t^2} -a^2\bigg{(}\sum_{i=1}^{d}\frac{1}{s^2_i}\frac{\partial^2 \tilde{u}_{NN}(\tilde{\bm{x}}_R^i,t_R^i)}{\partial \tilde{x}_i^2} \bigg{)} - \tilde{f}\big{(}\tilde{\bm{x}}_R^i,t_R^i,\tilde{u}_{NN}(\tilde{\bm{x}}_R^i,t_R^i)\big{)}\right|^2\\
    &+ \frac{1}{N_I}\sum_{i=1}^{N_I}\bigg{|}\tilde{u}_{NN}(\tilde{\bm{x}}_I^i, t_0) - \tilde{h}_K(\tilde{\bm{x}}_I^i)\bigg{|}^2 + \frac{1}{N_I}\sum_{i=1}^{N_I}\bigg{|}\frac{\partial \tilde{u}_{NN}(\tilde{\bm{x}}_I^i, t_0)}{\partial t} - \tilde{h}_P(\tilde{\bm{x}}_I^i)\bigg{|}^2\\
    &+ \frac{1}{N_B}\sum_{i=1}^{N_B}\bigg{|}\tilde{u}_{NN}\left(\tilde{\bm{x}}_B^i,t_B^i\right)-\tilde{g}_D(\tilde{\bm{x}}_B^i,t_B^i)\bigg{|}^2
\end{aligned}
\end{equation}
And the loss function of S-NFPINN method for the wave problem with Neumann boundary is formulated according to Eqs. \eqref{normalized:wave_space}, \eqref{normalized:wave_space:Neumman} and \eqref{normalized:wave_space:initial}, it is as follows
\begin{equation}\label{Loss_NFPINN_Neumman}
\begin{aligned}
    Loss_{Neu}&=\frac{1}{N_R}\sum_{i=1}^{N_R}\left| \frac{\partial^2 \tilde{u}_{NN}(\tilde{\bm{x}}_R^i,t_R^i)}{\partial t^2} -a^2\bigg{(}\sum_{i=1}^{d}\frac{1}{s^2_i}\frac{\partial^2 \tilde{u}_{NN}(\tilde{\bm{x}}_R^i,t_R^i)}{\partial \tilde{x}_i^2} \bigg{)} - \tilde{f}\big{(}\tilde{\bm{x}}_R^i,t_R^i,\tilde{u}_{NN}(\tilde{\bm{x}}_R^i,t_R^i)\big{)}\right|^2\\
    &+ \frac{1}{N_I}\sum_{i=1}^{N_I}\bigg{|}\tilde{u}_{NN}(\tilde{\bm{x}}_I^i, t_0) - \tilde{h}_K(\tilde{\bm{x}}_I^i)\bigg{|}^2 + \frac{1}{N_I}\sum_{i=1}^{N_I}\bigg{|}\frac{\partial \tilde{u}_{NN}(\tilde{\bm{x}}_I^i, t_0)}{\partial t} - \tilde{h}_P(\tilde{\bm{x}}_I^i)\bigg{|}^2\\
    &+ \frac{1}{N_B}\sum_{i=1}^{N_B}\bigg{|}\frac{\partial \tilde{u}_{NN}(\tilde{\bm{x}}_B^i,t_B^i)}{\partial \vec{n}}-\tilde{g}_N(\tilde{\bm{x}}_B^i,t_B^i)\bigg{|}^2
\end{aligned}
\end{equation}

Similarly, the loss function of NFPINN method for Dirichlet and Neumman boundaries can be obtained according to the normalized formulation of wave problems that adopt temporal normalization technique and temporal-spatial normalization technique, respectively.

Therefore, our objective is to determine an optimal set of parameters, denoted as $\bm{\theta}^*$, such that the approximations $\tilde{u}_{NN}$ effectively minimize the loss function $Loss(\bm{x};\bm{\theta})$. To achieve this ideal set of parameters, one can iteratively refine the weights and biases of DNN by employing optimization techniques, including gradient descent (GD) or stochastic gradient descent (SGD), throughout the training process. Within this framework, the SGD method, which utilizes a "mini-batch" of training data to enhance computational efficiency and stability, is formulated as follows:

\begin{equation}\label{optimize}
\bm{\theta}^{k+1}=\bm{\theta}^{k}-\alpha^k\nabla_{\bm{\theta}^k}Loss(\bm{x}^k;\bm{\theta}^{k})~~\text{with}~~\bm{x}^k\in \mathcal{S}^k,
\end{equation}
where the ``learning rate'' $\alpha^k$ decreases with $k$ increasing, and $Loss(\bm{x}^k;\bm{\theta}^{k})=Loss_{Dir}$ or $Loss(\bm{x}^k;\bm{\theta}^{k})=Loss_{Neu}$. 

\begin{remark}
The training loss of the NFPINN method for the $k_{th}$ iteration is defined as
$\varepsilon_T(\bm{\theta},\mathcal{S}^k)$, it describes the mapping $\bm{\theta}\mapsto\varepsilon_T(\bm{\theta},\mathcal{S}^k)$, then we can get
\begin{equation*}
    \varepsilon_T(\bm{\theta},\mathcal{S}^k) = Loss_{Dir} ~~\text{or}~~\varepsilon_T(\bm{\theta},\mathcal{S}^k) = Loss_{Neu}
\end{equation*}

Furthermore, the generalization error associated with our spatially normalized method, when applied to the wave problem \eqref{eq:wave equation} with Dirichlet boundary conditions, can be expressed in the following form.

\begin{equation}\label{gen_err}
    \begin{aligned}
     \varepsilon_G(\bm{\theta})&=\int\int_{\overline{\Omega}\times T}\bigg{|}\frac{\partial^2 \tilde{u}_{NN}(\tilde{\bm{x}}, t)}{\partial t^2} - a^2\bigg{(}\sum_{i=1}^{d}\frac{1}{s^2_i}\frac{\partial^2 \tilde{u}_{NN}(\tilde{\bm{x}}, t)}{\partial \tilde{x}_i^2} \bigg{)} - \tilde{f}\big{(}\tilde{\bm{x}}, t, \tilde{u}_{NN}(\tilde{\bm{x}}, t)\big{)}\bigg{|}^2d\bm{x}dt\\
     &+\int\int_{\partial\overline{\Omega}\times T}\bigg{|}\tilde{u}_{NN}(\tilde{\bm{x}},t) - \tilde{g}_D(\tilde{\bm{x}},t)\bigg{|}^2d\bm{x}dt\\
     &+\int_{\overline{\Omega}}\bigg{|}    \tilde{u}_{NN}(\tilde{\bm{x}},t_0) - \tilde{h}_K(\tilde{\bm{x}})\bigg{|}^2d\bm{x} +\int_{\overline{\Omega}}\bigg{|}\frac{\partial \tilde{u}_{NN}(\tilde{\bm{x}},t_0)}{\partial t} - \tilde{h}_P(\tilde{\bm{x}})\bigg{|}^2d\bm{x}
    \end{aligned}
\end{equation}

Similarly, the general error for other cases can be derived based on the corresponding equations for T-NFPINN and ST-NFPINN, taking into account whether the boundary constraints are of the Dirichlet or Neumann type.

The total error of the NFPINN approximation is defined in terms of $\hat{u} = \tilde{u}_{NN}-\tilde{u}$, which represents the discrepancy between the exact solution of the wave equation and the corresponding approximation obtained through the NFPINN method.

\begin{equation}\label{err2pinn}
 \varepsilon(\bm{\theta}) = \int_{\Omega}\bigg{(}|\hat{u}(\bm{x}, t)|+|\nabla \hat{u}(\bm{x}, t)|\bigg{)}d\bm{x}dt
\end{equation}

According to the analysis presented in \cite{qian2023physics}, the total error $\varepsilon(\bm{\theta})^2$ remains small when the generalization error $\varepsilon_G(\bm{\theta})^2$ is sufficiently small for the NFPINN approximation $\tilde{u}(\cdot;\bm{\theta})$ . Furthermore, it can be demonstrated that the total error 
$\varepsilon(\bm{\theta})^2$ can be made extremely small, provided that the training error $\varepsilon_T(\bm{\theta},\mathcal{S})^2$ is minimized to a sufficiently small value and the sample set is sufficiently large. Consequently, the total error is inherently bounded, ensuring that our method yields reliable and satisfactory results for solving the wave equation.
\end{remark}

\section{The process of NFPINN algorithm}\label{sec:4}
For the NFPINN method, where the solver is an FFM-based DNN (FFM-DNN) consisting of $Q$ subnetworks, as illustrated in Fig.\ref{Fig:FPINN network}, the input data for each subnetwork undergoes a translation through the following operation.
\begin{equation}
\hat{\bm{x}}=a_i*\bm{x}~~i=1, 2, \dots, Q
\end{equation}
where $a_i\geq 1$ represents a positive scalar factor, implying that the scale vector is defined as $\Lambda = (a_1, a_2,\cdots, a_Q)$ and denoting the output of each subnetwork as $\bm{F}_i(i=1, 2, \cdots, Q)$, the overall output of the FFM-DNN model is consequently derived as
\begin{equation}
\bm{y}(\bm{x};\bm{\theta})=\frac{1}{Q}\sum_{i=1}^{Q}\bm{F}_i
\end{equation}

Based on the aforementioned discussions, the procedure of the NFPINN algorithm, which is designed to address the wave equation within finite-dimensional spaces, is outlined in detail as follows.

\begin{algorithm}[H]
    \caption{The NFPINN algorithm designed for addressing the wave equation}
    1. Generating the $k_{th}$ training set $\mathcal{S}^k=\mathcal{S}_R^k\cup\mathcal{S}_B^k\cup\mathcal{S}_I^k$ includes interior points $S_{R}^k=\{(\bm{x}^i_R,t^i_R)\}_{i=1}^{N_{R}}$ with $\bm{x}_R^i\in\mathbb{R}^d$ and $t_R^i\in\mathbb{R}$, boundary points $S_{B}^k=\{(\bm{x}^j_B,t^j_B)\}_{j=1}^{N_{B}}$ with $\bm{x}_B^j\in\mathbb{R}^d$ and $t_B^j\in\mathbb{R}$, initial points $S_{I}^k=\{(\bm{x}^j_I,t_0)\}_{j=1}^{N_{B}}$ with $\bm{x}_I^j\in\mathbb{R}^d$. Here, we draw the above points by three method with positive probability density $\nu$, such as Latin hypercube sampling approach.
    
    2. Obtaining the objective function $Loss(\mathcal{S}^k;\bm{\theta}^{k})$ for on the training set $\mathcal{S}^k$:
    \begin{equation*}
    Loss(\mathcal{S}^k;\bm{\theta}^{k}) =Loss_{Dir}~~\text{or}~~Loss(\mathcal{S}^k;\bm{\theta}^{k}) =Loss_{Neu}
    \end{equation*}
    with $Loss_{Dir}$ and $Loss_{Neu}$ being defined in Eq.\eqref{Loss_NFPINN_Dirichlet} and Eq.\eqref{Loss_NFPINN_Neumman}, respectively.
    
    3. Adopting an appropriate optimization method to adjust the internal parameters of FFM-DNN at the point $(\tilde{\bm{x}}^k, \tilde{t}^k)$ in $\mathcal{S}^k$. One possible choice is the SGD method described as follows.
    
    \begin{equation*}
    \bm{\theta}^{k+1}=\bm{\theta}^{k}-\alpha^k\nabla_{\bm{\theta}^k}Loss(\tilde{\bm{x}}^k,\tilde{t}^k;\bm{\theta}^{k}),
    \end{equation*}
    where the ``learning rate'' $\alpha^k$ reduces as $k$ increasing.
    
    4. Steps 1-3 are iteratively executed in succession until the convergence criterion is met or the objective function exhibits a tendency toward stabilization.
    \label{algor:FMPINN}
\end{algorithm}

\section{Numerical experiments}\label{sec:5}
In this section, we conduct a comprehensive evaluation of the performance of the NFPINN method by implementing three distinct normalization strategies to solve the wave propagation equation, considering various boundary conditions across two-dimensional and three-dimensional spaces.

\subsection{Training settings}
\begin{itemize}
    \item 
    \emph{Network configuration}: To ensure a fair comparison, all four models discussed in the previous sections share a consistent network configuration, wherein the solvers for the FPINN and NFPINN methods are set as a 4-layer FFM-based DNN  with layer sizes of (20, 15, 15, 10), employing the GELU activation function across all hidden layers while utilizing a linear output layer for the final predictions. The training process is conducted using the Adam optimizer, which operates under an exponential learning rate strategy, starting with an initial learning rate of 0.01 that undergoes a decay of 0.035 every 100 training epochs, ultimately spanning a total of 30,000 epochs for all models.
    \item 
    \emph{Sampling strategy}: During each training epoch, Latin Hypercube Sampling (LHS) is employed as a strategic method to intelligently select training points, ensuring a well-distributed representation across both the non-normalized and normalized domains.
    \item 
    \emph{Criterion selection}: To evaluate the accuracy of the aforementioned models, we utilize the Relative Error (REL) as the evaluation metric, which is defined as follows.
    \begin{equation}
    REL = \sqrt{\frac{\sum_{i=1}^{N'}\bigg{(}u_{NN}(x^i,t^i)-u(x^i,t^i)\bigg{)}^2}{\sum_{i=1}^{N'}\bigg{(}u(x^i,t^i)\bigg{)}^2}}
    \end{equation}
    where $u_{NN}(x^i,t^i)$ denotes the approximate solution predicted by the neural network, and $u(x^i,t^i)$ corresponds to the exact or reference solution, with the set ${(x^i,t^i)}^{N'}_{i=1}$ comprising the testing points. Here, the parameter \emph{N'} signifies the total number of testing points.
    \item 
    \emph{Visualization}: The REL of the test set is recorded at intervals of 1,000 steps, with its trends visually illustrated through graphical representation. Furthermore, the point-wise absolute errors between the prediction outcomes of the FPINN and NFPINN models are also depicted graphically to provide a comparative analysis.
    \item 
    \emph{Computing power resources}: Implemented using PyTorch, our framework runs on a desktop computer equipped with 16 GB of RAM and a single NVIDIA GeForce GTX 1080 GPU with 8 GB of memory, providing the necessary computational resources to ensure both the efficient and accurate training of the FPINN and NFPINN models while enabling a comprehensive analysis of their performance.
\end{itemize}

\subsection{Numerical case studies and corresponding outcomes}
In this section, we compare the performance of FPINN, S-NFPINN, T-NFPINN, ST-NFPINN models for 2D and 3D wave equations with Dirichlet and Neumann boundary conditions.
\begin{example}
We aim to approximate the solution of Eq.\eqref{eq:Example_1_govern} with Dirichlet boundary for given regular domain $\Omega=[0,10\pi]\times[0,10\pi], t\in(0,10)$. 
\end{example}

\begin{figure}[H]
    \centering
    \subfigure[Non-normalized domain]{
        \label{E1-09} 
        \includegraphics[scale=0.325]{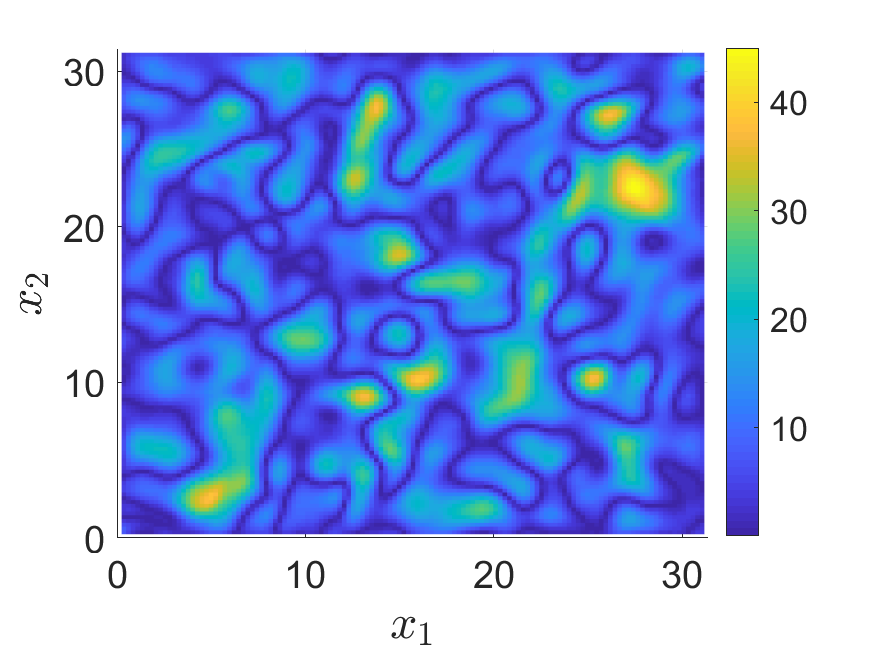}
    }
    \subfigure[Normalized space domain]{
        \label{E1-10} 
        \includegraphics[scale=0.325]{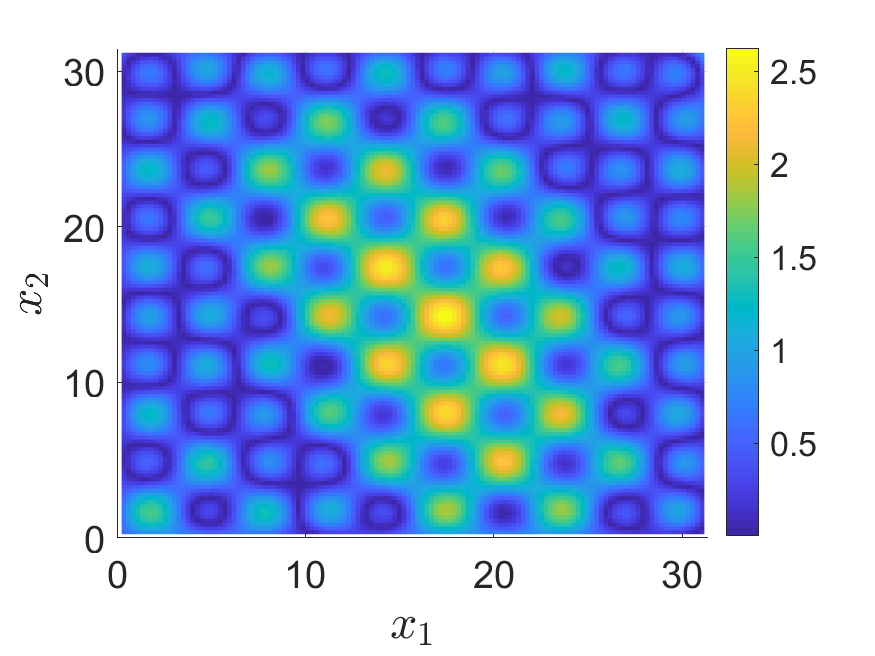}
    }
    \subfigure[Normalized time domain]{
        \label{E1-11} 
        \includegraphics[scale=0.325]{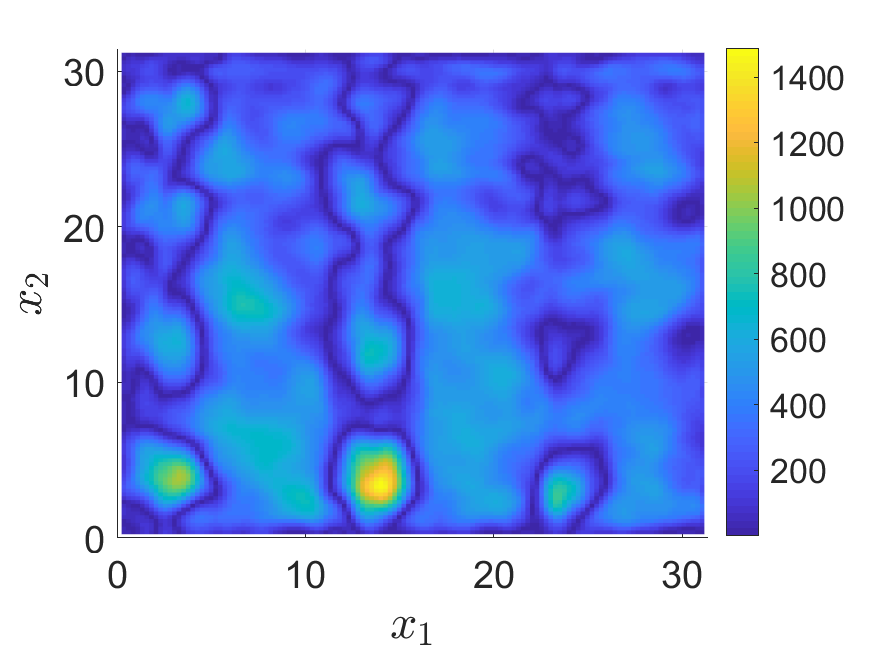}
    }
    \subfigure[Normalized space and time domain]{
        \label{E1-12}
        \includegraphics[scale=0.325]{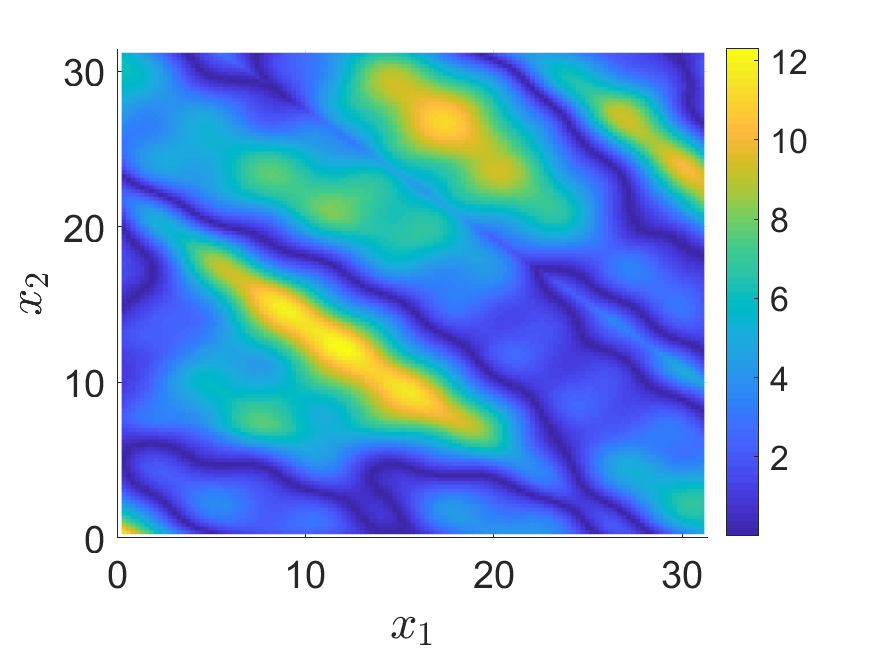}
    }
    \subfigure[Relative error]{
        \label{E1-13}
        \includegraphics[scale=0.31]{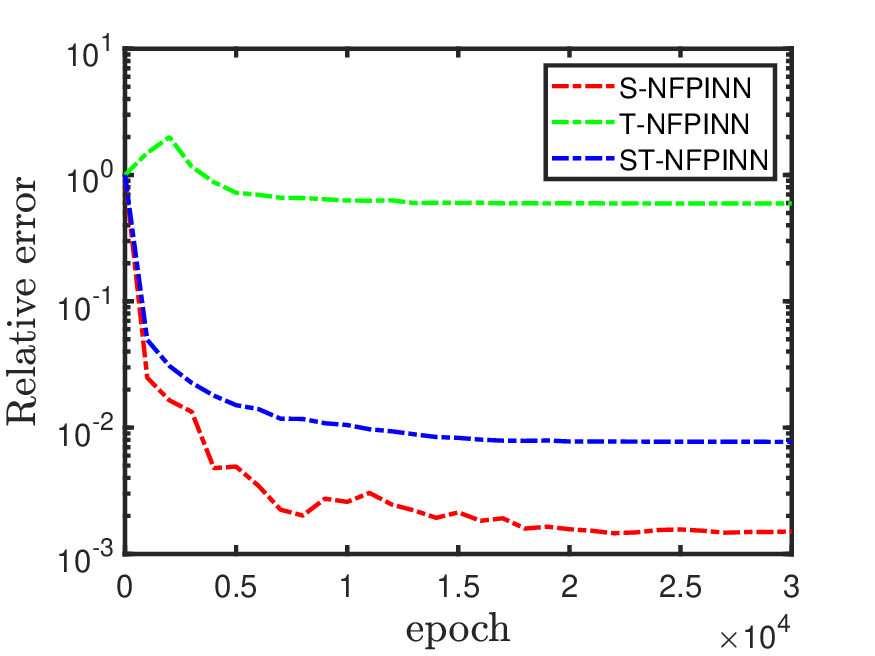}
    }    
    \caption{Numerical results of NFPINN in different normalized methods for test example in domain $\Omega$.}
    \label{Fig:Example_1_NFPINN}
\end{figure}

We approximate the solution of Eq.\eqref{eq:Example_1_govern} using the three NFPINN methods on 16384 equidistant grid points distributed across regular domain. In the training set, we have carefully chosen 700 data points for each spacetime under initial conditions, along with 1,500 points distributed throughout the domain and 300 points specifically allocated for boundary conditions. Fig.\ref{Fig:Example_1_NFPINN} illustrates the point-wise absolute errors distribution of the three NFPINN models across the domain. 

The reduction in color intensity of the point-wise absolute errors heatmaps for the four models is demonstrated in Fig.\ref{E1-09} to Fig.\ref{E1-12}. After training using the S-NFPINN method, the absolute errors of the model drop a lot, as compared to the FPINN method. The relative error curves of S-NFPINN, T-NFPINN and ST-NFPINN methods are shown in Fig.\ref{E1-13}. The S-NFPINN model's accuracy jumps to $10^{-3}$ and exhibits a more rapid convergence rate.. This is consistent with our study that utilizing a normalized technology during training can be more effective in addressing the non-unitized domain and long time range. To sum up, S-NFPINN outperforms FPINN, T-NFPINN and ST-NFPINN under the two-dimensional wave equation with Dirichlet boundary conditions.

\begin{example}
In order to make the normalized method more convincing, we solve a regular high-frequency 2D wave equation in $\Omega=[0,10\pi]\times[0,10\pi], t\in(0,10)$. The govern conditions is described as
\begin{equation}\label{eq:Example_2_govern}
    \frac{\partial^2 u}{\partial t^2} =  0.01*\bigg{(}\frac{\partial^2 u}{\partial x_1^2} + \frac{\partial^2 u}{\partial x_2^2}\bigg{)} + 12t^2
\end{equation}

Boundary conditions are
\begin{equation}
    \begin{aligned}
        u|_{x_1=0}(x_1,x_2,t) &= u|_{x_1=10\pi}(x_1,x_2,t) = t^4 + \cos(6*x_2)\cdot\sin(t)\\
        u|_{x_2=0}(x_1,x_2,t) &= u|_{x_2=10\pi}(x_1,x_2,t) = t^4 + \cos(8*x_1)\cdot\sin(t)
    \end{aligned}
\end{equation}

Initial conditions are
\begin{equation}
    \begin{aligned}
        u|_{t=0}(x_1,x_2,t) &= 0\\
        \frac{\partial u}{\partial t}\bigg{|}_{t=0}(x_1,x_2,t) &= \cos(8*x_1)\cdot\cos(6*x_2)
    \end{aligned}
\end{equation}

The analytical solution is
\begin{equation}
    u(x_1,x_2,t) = t^4 + \cos(8*x_1)\cdot\cos(6*x_2)\cdot\sin(t)
\end{equation}
\end{example}

By selecting 700 points for each spacetime under initial conditions, 1500 points within the domain, and 300 points under boundary conditions, the training set is meticulously constructed to ensure that the models are not only well-trained but also capable of accurately capturing intricate spatiotemporal patterns. The four model is composed of 10 subnetworks according to the manually defined frequencies $\bm{\Lambda}=(1,2,3,4,5,6,7,8,9,10)$. Additionally, to rigorously assess the feasibility of the model, a testing set comprising 16,384 sampled collocation points is drawn from the regular domain $\Omega$ at $t = 2.5$. This strategic selection ensures a comprehensive evaluation of the model's generalization capability and its ability to accurately capture the underlying spatiotemporal dynamics.

\begin{figure}%[H]
    \centering
    \subfigure[Analytical solution]{
        \label{E2-01} 
        \includegraphics[scale=0.325]{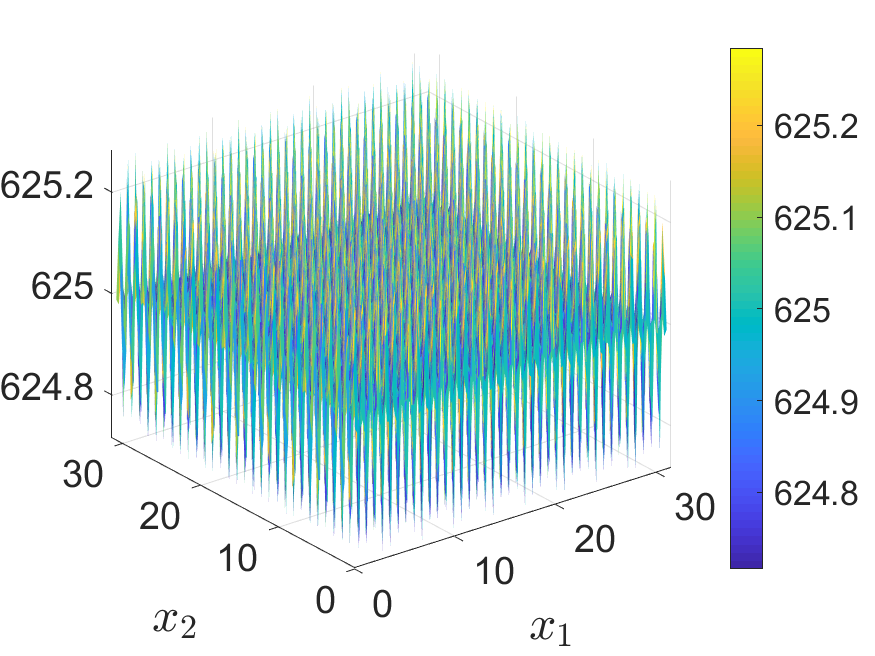}
    }
    \subfigure[Normalized space domain]{
        \label{E2-02} 
        \includegraphics[scale=0.325]{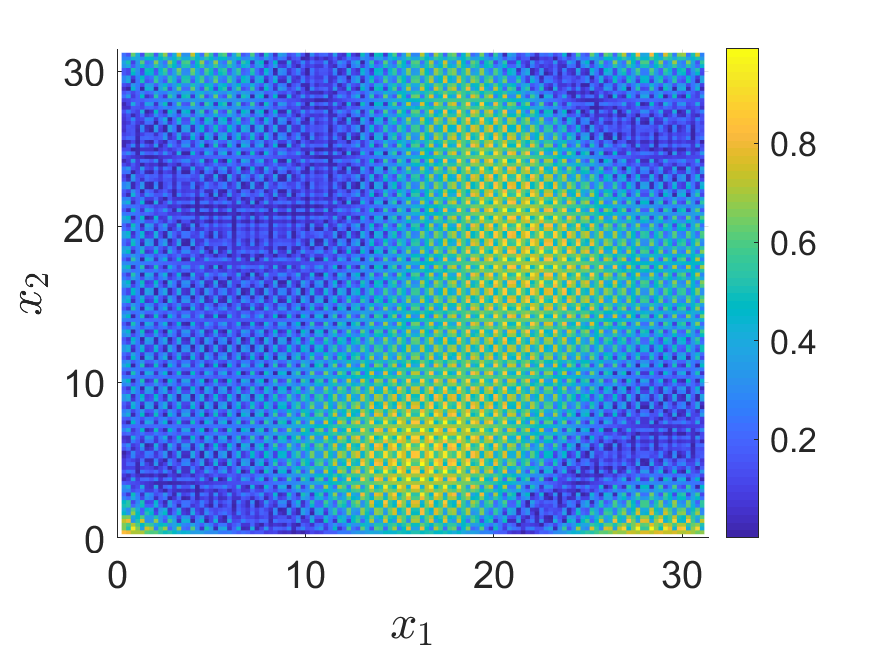}
    }
    \subfigure[Normalized time domain]{
        \label{E2-03}
        \includegraphics[scale=0.325]{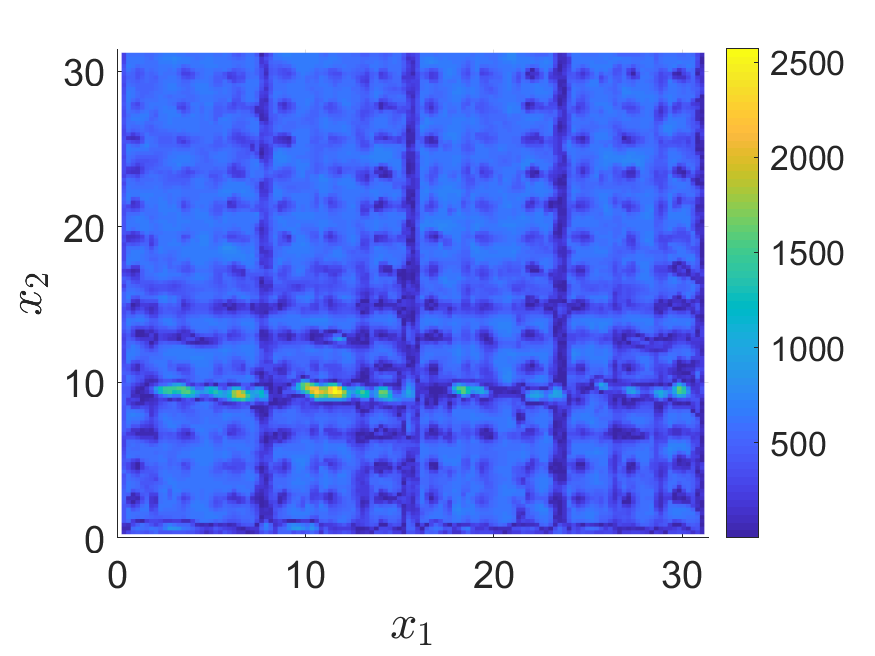}
    }
    \subfigure[Normalized space and time domain]{
        \label{E2-04}
        \includegraphics[scale=0.325]{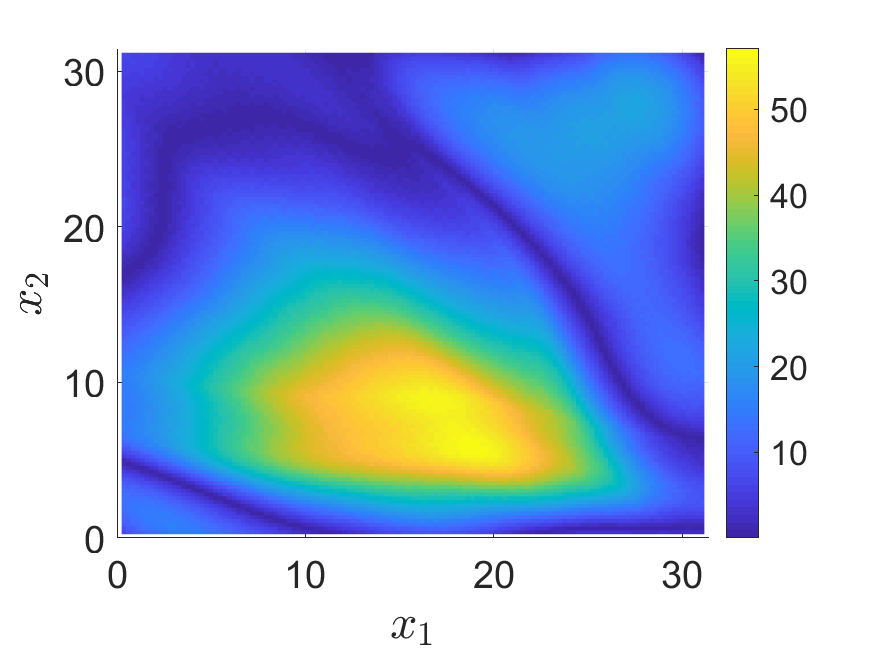}
    }  
    \subfigure[Relative error]{
        \label{E2-05}
        \includegraphics[scale=0.305]{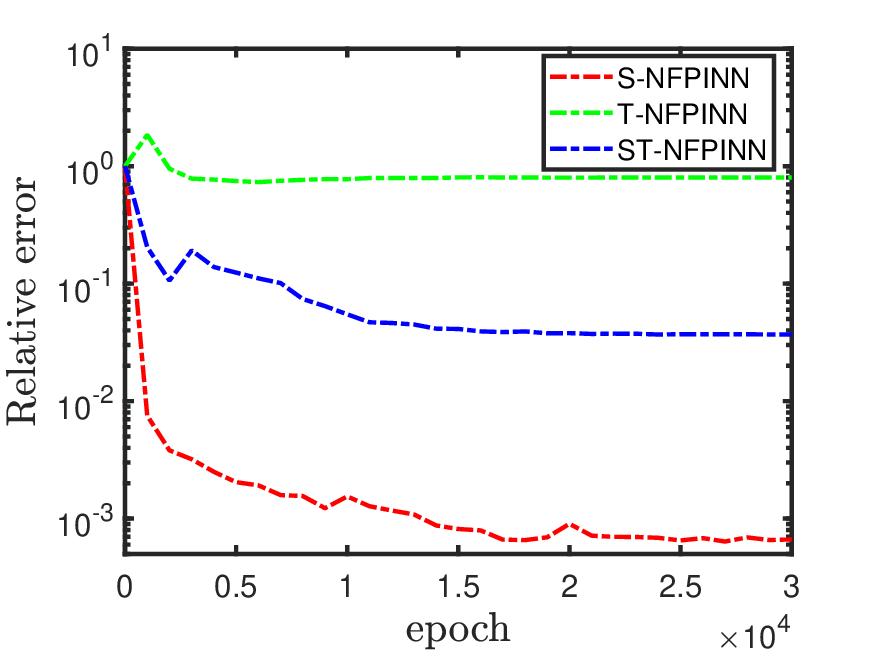}
    }    
    \caption{Numerical results of NFPINN method for test example 2.}
    \label{Fig:Example_2_NFPINN}
\end{figure}

The numerical results shown in Fig.\ref{E2-01}-Fig.\ref{E2-05} provide compelling evidence that S-NFPINN remains the optimal choice. This conclusion is drawn from its outstanding performance in terms of both point-wise absolute errors and relative errors. The data clearly illustrate that S-NFPINN consistently surpasses T-NFPINN and ST-NFPINN methods, delivering greater accuracy and reliability across multiple test cases. These results emphasize the robustness and efficiency of S-NFPINN.

\begin{example}
Consider a porous 2D wave equation with Neumann boundary condition in $\Omega=[0,10\pi]\times[0,10\pi], t\in(0,10)$. The governed equation is
\begin{equation}\label{eq:Example_3_govern}
    \frac{\partial^2 u}{\partial t^2} = \frac{1}{2}\bigg{(}\frac{\partial^2 u}{\partial x_1^2} + \frac{\partial^2 u}{\partial x_2^2}\bigg{)} + 12t^2 - \cos(x_1)\cdot\sin(x_2)\cdot\sin(t)
\end{equation}

Boundary conditions are
\begin{equation}
    \begin{aligned}
        \frac{\partial u}{\partial x_1}\bigg{|}_{x_1=0}(x_1,x_2,t) &= 0 \\
        \frac{\partial u}{\partial x_1}\bigg{|}_{x_1=10\pi}(x_1,x_2,t) &= \sin(x_2)\cdot\sin(t)\cdot{10\pi}  \\
        \frac{\partial u}{\partial x_2}\bigg{|}_{x_2=0}(x_1,x_2,t) &= \frac{\partial u}{\partial x_2}\bigg{|}_{x_2=10\pi}(x_1,x_2,t) =\sin(x_1)\cdot\sin(t)\cdot{x_1} \\ 
    \end{aligned}
\end{equation}

Initial conditions are
\begin{equation}
    \begin{aligned}
        u|_{t=0}(x_1,x_2,t) &= 0  \\
        \frac{\partial u}{\partial t}\bigg{|}_{t=0}(x_1,x_2,t) &= \sin(x_1)\cdot\sin(x_2)\cdot{x_1}
    \end{aligned}
\end{equation}

The analytical solution is
\begin{equation}
    u(x_1,x_2,t) = t^4 + \sin(x_1)\cdot\sin(x_2)\cdot\sin(t)\cdot{x_1}
\end{equation}
\end{example}

Within the training set, we have meticulously selected 700 data points for each spacetime under initial conditions, 1500 points distributed across the domain, and 300 points designated for boundary conditions. This careful selection guarantees that the models undergo comprehensive training, enabling them to effectively capture intricate spatiotemporal dynamics. Furthermore, the four models are structured as an ensemble of 10 interconnected subnetworks, each specifically adjusted to the manually defined set of frequencies $\bm{\Lambda}=(1,2,3,4,5,6,7,8,9,10)$, thereby enhancing their capacity to learn and represent complex patterns. Moreover, to thoroughly assess the feasibility of the model, a testing set consisting of 17,627 data points, strategically sampled from the porous domain $\Omega$ at $t = 2.5$, is utilized. This carefully designed evaluation framework ensures a rigorous validation of the model's ability to generalize and accurately capture the underlying spatiotemporal characteristics of the system.

\begin{figure}%[H]
    \centering
    \subfigure[Analytical solution]{
        \label{E3-01} 
        \includegraphics[scale=0.325]{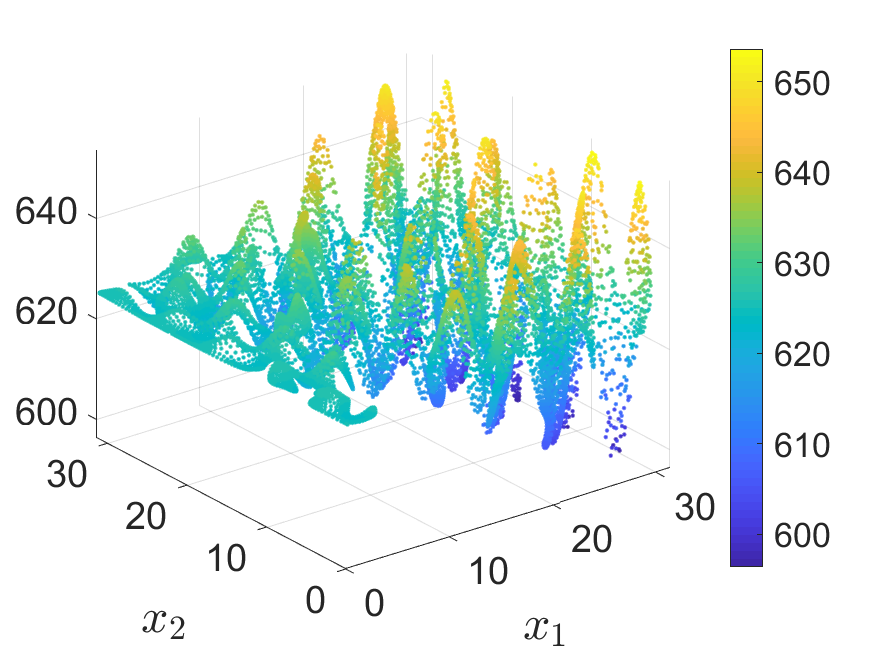}
    }
    \subfigure[Normalized space domain]{
        \label{E3-02} 
        \includegraphics[scale=0.325]{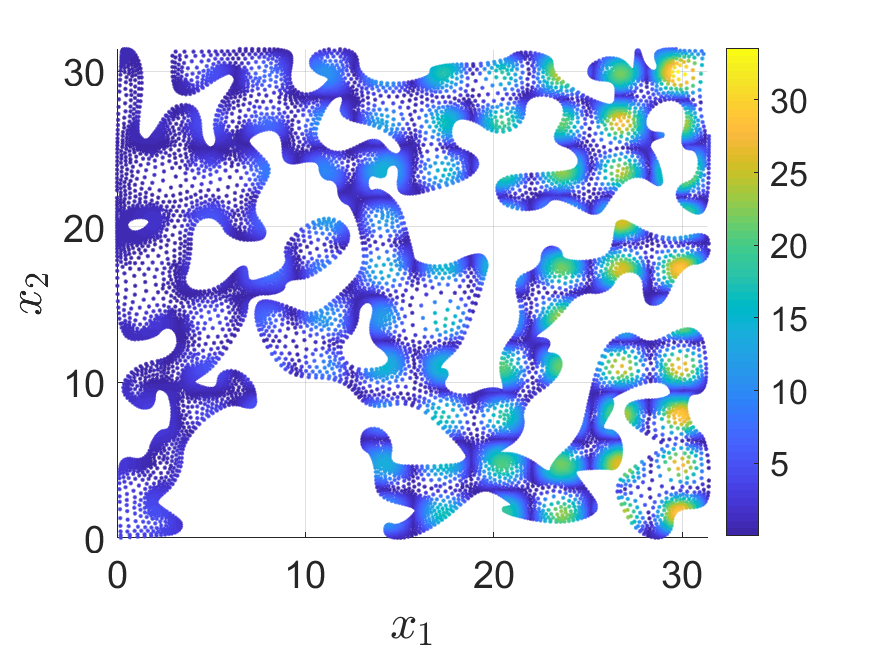}
    }
    \subfigure[Normalized time domain]{
        \label{E3-03}
        \includegraphics[scale=0.325]{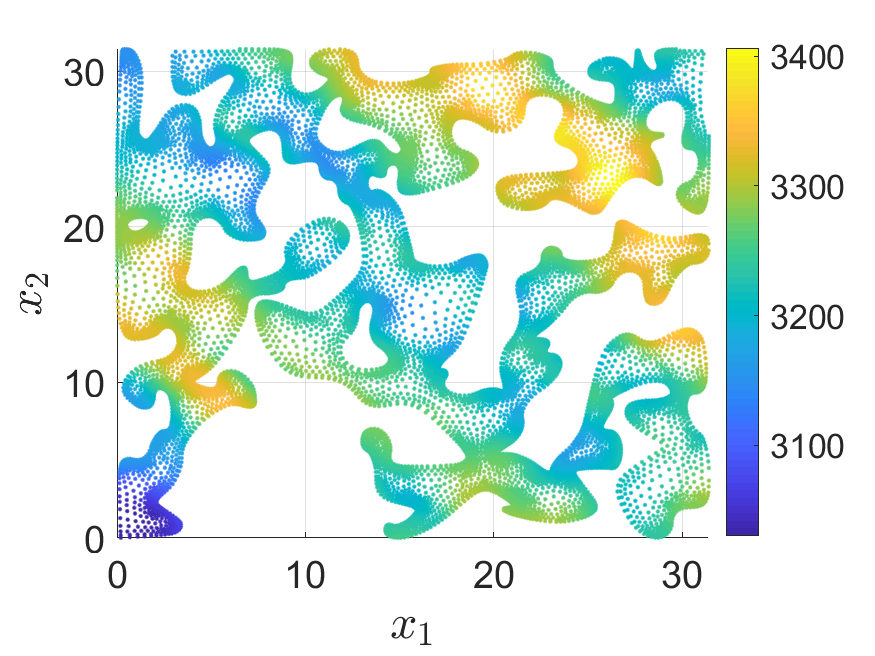}
    }
    \subfigure[Normalized space and time domain]{
        \label{E3-04}
        \includegraphics[scale=0.325]{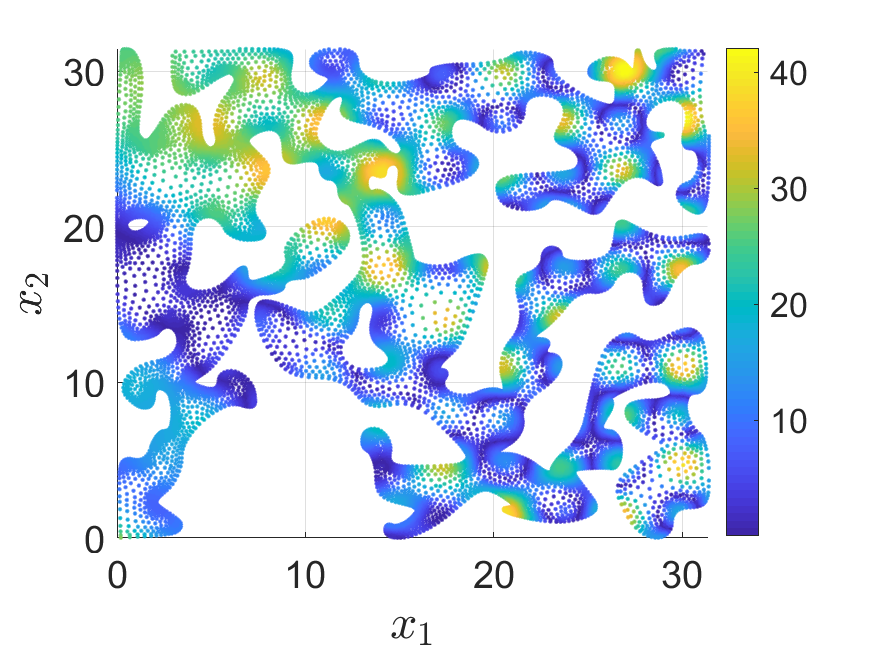}
    }  
    \subfigure[Relative error]{
        \label{E3-05}
        \includegraphics[scale=0.305]{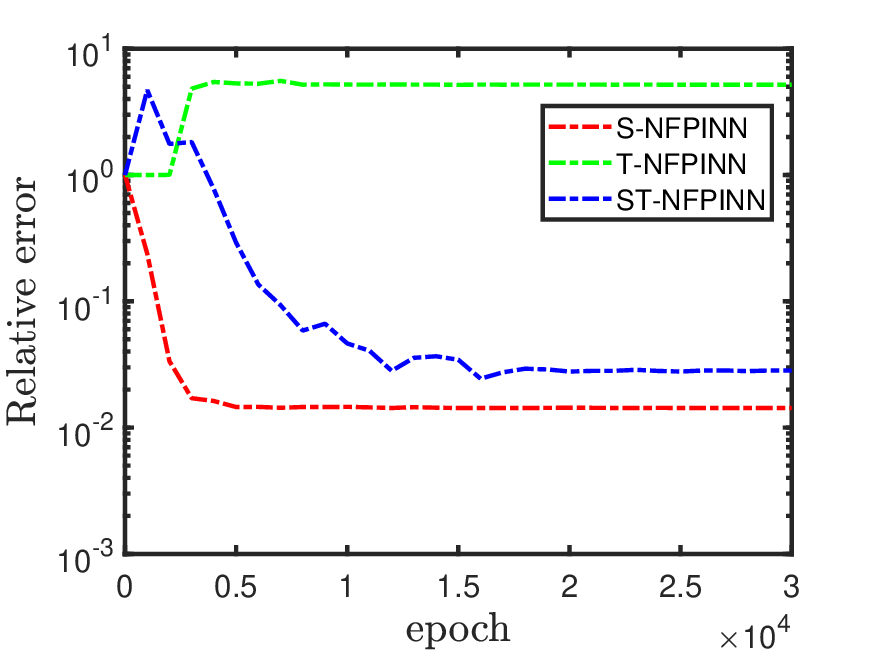}
    }    
    \caption{Numerical results of NFPINN method for test example 3.}
    \label{Fig:Example_3_NFPINN}
\end{figure}

In Fig.\ref{E3-01}–Fig.\ref{E3-04}, we present the discrepancies between the approximate solutions and the exact solutions, illustrating the accuracy of S-NFPINN method. Additionally, the relative errors are depicted in Fig.\ref{E3-05}, providing a quantitative assessment of the approximation quality. The results in Fig.\ref{Fig:Example_3_NFPINN} highlight the superior performance of the S-NFPINN method in solving Eq.\eqref{eq:Example_3_govern}, particularly in handling complex boundary conditions within an irregular two-dimensional domain. Compared to T-NFPINN and ST-NFPINN, S-NFPINN demonstrates a significantly improved capability in capturing the solution with higher precision and robustness. This advantage is attributed to its enhanced ability to adapt to intricate geometries and enforce boundary constraints more effectively.  The reduced error and better convergence behavior further confirm the efficacy of S-NFPINN in solving complex partial differential equations.

\begin{example}
Now, we extend from 2D to 3D to solve a spherical cavity wave equation under the Dirichlet boundary condition. The problem domain is depicted as $\Omega=[0,10\pi]\times[0,10\pi]\times[0,10\pi], t\in(0,10)$, respectively. The governed equation is

\begin{equation}\label{eq:Example_4_govern}
    \frac{\partial^2 u}{\partial t^2} = \frac{1}{2}\bigg{(}\frac{\partial^2 u}{\partial x_1^2} + \frac{\partial^2 u}{\partial x_2^2} + \frac{\partial^2 u}{\partial x_3^2}\bigg{)} + 12t^2 + \frac{1}{2}\cdot\sin(x_3) 
\end{equation}
    
Boundary conditions are
\begin{equation}
    \begin{aligned}
        u|_{x_1=0}(x_1,x_2,x_3,t) &= u|_{x_1=10\pi}(x_1,x_2,x_3,t) = t^4 + \sin(x_3)\\
        u|_{x_2=0}(x_1,x_2,x_3,t) &= u|_{x_2=10\pi}(x_1,x_2,x_3,t) = t^4 + \sin(x_3)\\
        u|_{x_3=0}(x_1,x_2,x_3,t) &= u|_{x_3=10\pi}(x_1,x_2,x_3,t) = t^4 + \sin(x_1)\cdot\sin(x_2)\cdot\sin(t)
    \end{aligned}
\end{equation}
    
Initial conditions are
\begin{equation}
    \begin{aligned}
        u|_{t=0}(x_1,x_2,x_3,t) &= \sin(x_3)\\
        \frac{\partial u}{\partial t}\bigg{|}_{t=0}(x_1,x_2,x_3,t) &= \sin(x_1) \cdot \sin(x_2)
    \end{aligned}
\end{equation}
    
The analytical solution is
\begin{equation}
    u(x_1,x_2,,x_3,t) = t^4 + \sin(x_1)\cdot\sin(x_2)\cdot\sin(t) + \sin(x_3)
\end{equation}
\end{example}

In constructing the training set, we have carefully selected 1000 data points for each spacetime under initial conditions, 2000 points distributed throughout the domain, and 500 points assigned to boundary conditions. This meticulous sampling strategy ensures that the models undergo thorough training, equipping them with the capability to accurately capture intricate spatiotemporal patterns. Moreover, the four models are structured as a composition of 15 interrelated subnetworks, each designed in accordance with the manually specified frequency set $\bm{\Lambda}=(1,2,3,4,5,6,7,8,9,10,11,12,13,14,15)$, thereby enhancing their ability to learn and represent complex dynamic behaviors. 

With the aim of ensuring a comprehensive evaluation of the model's performance in capturing the spatial and temporal characteristics of the domain, the testing set is constructed from 3,000 mesh grid points, which are systematically distributed across the spherical surface with a radius of $r=3$ at $t = 2.5$, thereby facilitating a detailed assessment of spatial and temporal variations.

\begin{figure}[H]
    \centering
    \subfigure[Analytical solution]{
        \label{E4-01} 
        \includegraphics[scale=0.325]{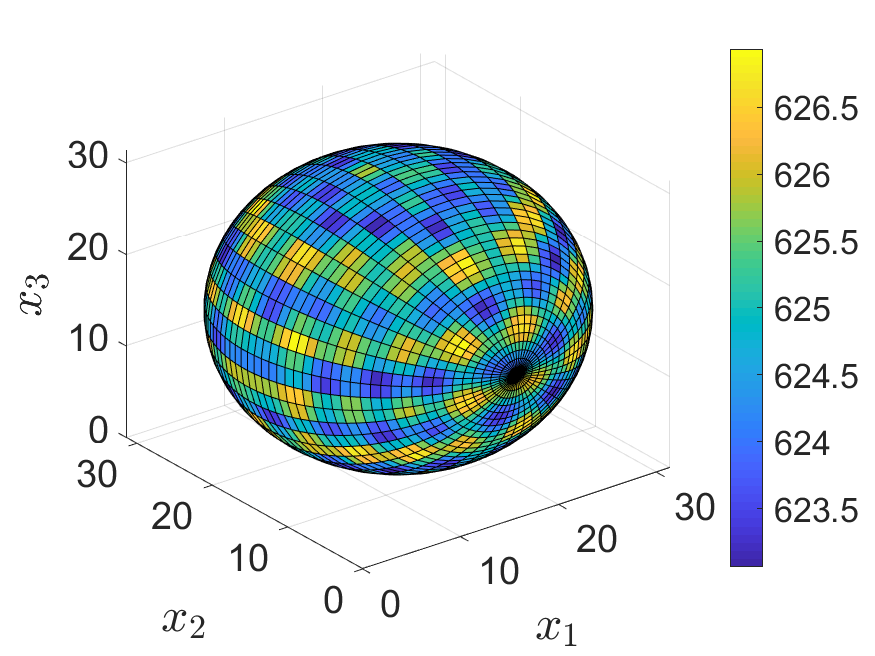}
    }
    \subfigure[Normalized space domain]{
        \label{E4-02} 
        \includegraphics[scale=0.325]{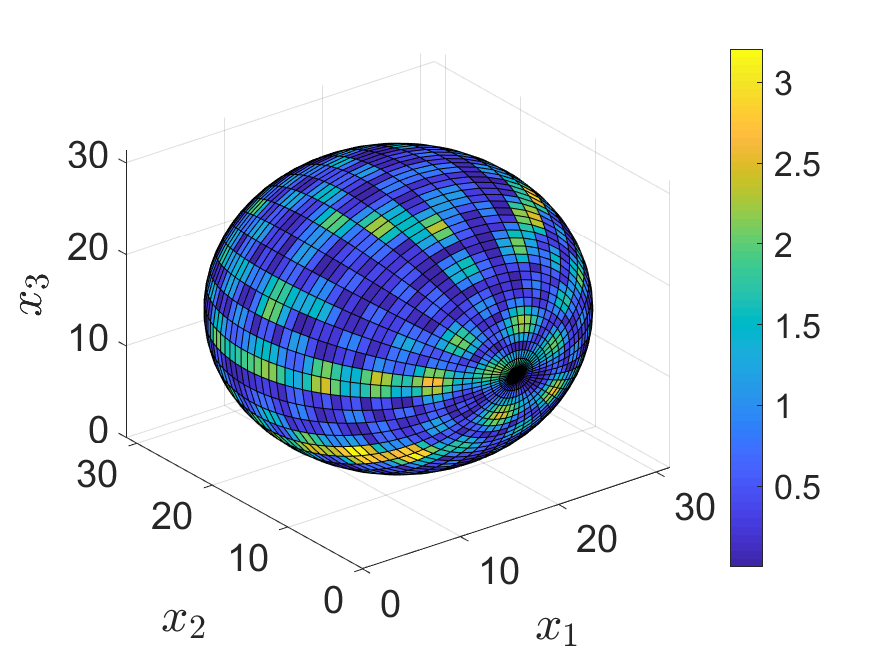}
    }
    \subfigure[Normalized time domain]{
        \label{E4-03}
        \includegraphics[scale=0.325]{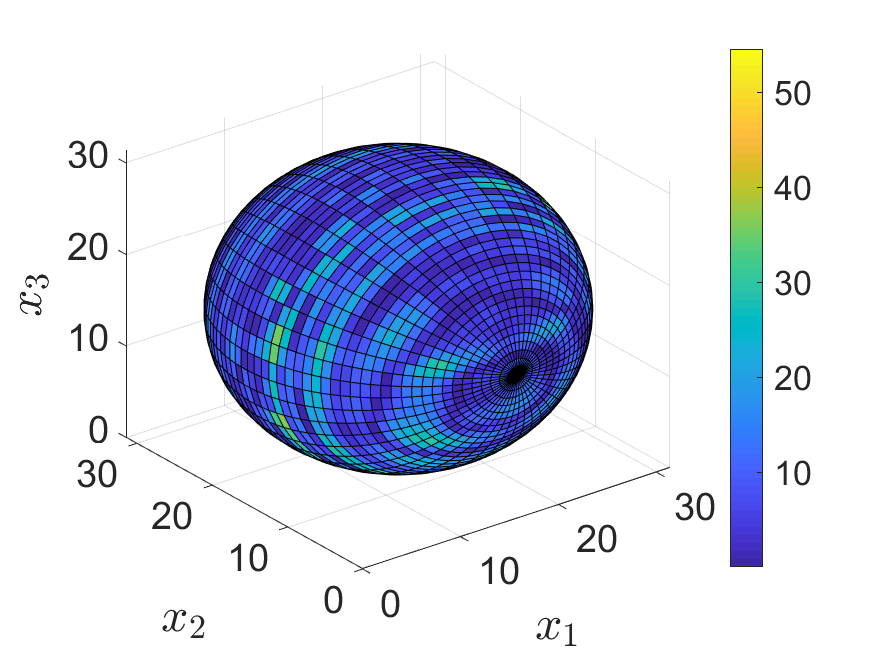}
    }
    \subfigure[Normalized space and time domain]{
        \label{E4-04}
        \includegraphics[scale=0.325]{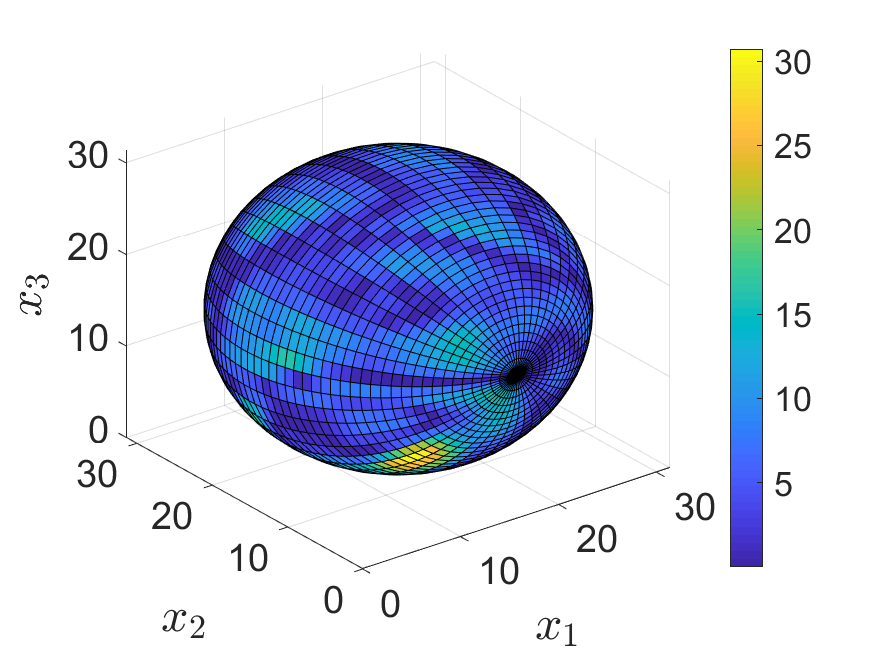}
    }  
    \subfigure[Relative error]{
        \label{E4-05}
        \includegraphics[scale=0.305]{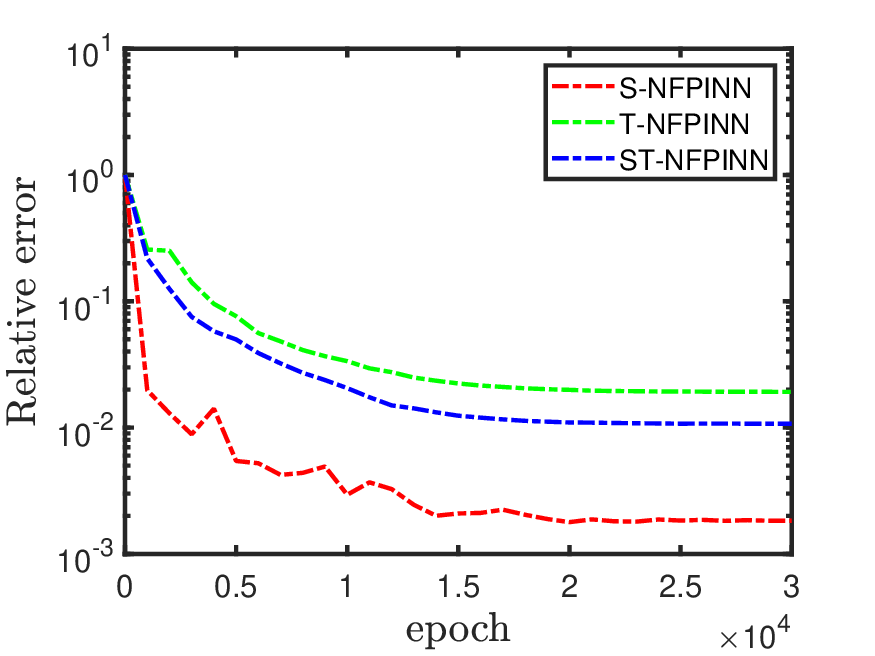}
    }    
    \caption{Numerical results of NFPINN method for test example 4.}
    \label{Fig:Example_4_NFPINN}
\end{figure}

For the spherical cavity in three-dimensional space, the point-wise errors presented in Fig.\ref{E4-01}- Fig.\ref{E4-04} clearly demonstrate that our proposed S-NFPINN method achieves significantly lower errors compared to T-NFPINN and ST-NFPINN methods when solving the wave equation Eq.\eqref{eq:Example_4_govern}. These results highlight the effectiveness of S-NFPINN in accurately capturing the wave propagation behavior while maintaining numerical stability across the entire computational domain. Notably, our method exhibits superior accuracy not only in the interior regions but also near the cavity boundary, where traditional approaches often struggle due to complex geometric constraints. This improvement is largely attributed to the integration of spectral neural functions with the physics-informed neural network framework, which enhances solution precision by enforcing physical constraints more effectively. Consequently, the results validate the robustness and efficiency of S-NFPINN as a powerful tool for solving wave equations in three-dimensional domains with intricate geometries.

\begin{example}
Similarly, we expand 2D to 3D to calculate porous wave equation with Neumann boundary condition in $\Omega=[0,10\pi]\times[0,10\pi]\times[0,10\pi], t\in(0,10)$. The governed equation is expressed as

\begin{equation}\label{eq:Example_5_govern}
    \frac{\partial^2 u}{\partial t^2} = \frac{\partial^2 u}{\partial x_1^2} + \frac{\partial^2 u}{\partial x_2^2} + \frac{\partial^2 u}{\partial x_3^2} + 12t^2 + 2\sin(x_1)\cdot\sin(x_2)
\end{equation}

The prescribe boundary conditions are
\begin{equation}
    \begin{aligned}
    \frac{\partial u}{\partial x_1}\bigg{|}_{x_1=0}(x_1,x_2,x_3,t) &= \frac{\partial u}{\partial x_1}\bigg{|}_{x_1=10\pi}(x_1,x_2,x_3,t) = \sin(x_2)\\
    \frac{\partial u}{\partial x_2}\bigg{|}_{x_2=0}(x_1,x_2,x_3,t) &= \frac{\partial u}{\partial x_2}\bigg{|}_{x_2=10\pi}(x_1,x_2,x_3,t) = \sin(x_1)\\
    \frac{\partial u}{\partial x_3}\bigg{|}_{x_3=0}(x_1,x_2,x_3,t) &= \frac{\partial u}{\partial x_3}\bigg{|}_{x_3=10\pi}(x_1,x_2,x_3,t) = \sin(t)  
    \end{aligned}
\end{equation}

and initial conditions are
\begin{equation}
    \begin{aligned}
        u|_{t=0}(x_1,x_2,x_3,t) &= \sin(x_1)\cdot\sin(x_2)\\
        \frac{\partial u}{\partial t}\bigg{|}_{t=0}(x_1,x_2,x_3,t) &= \sin(x_3)
    \end{aligned}
\end{equation}

respectively. An analytical solution is given by
\begin{equation}
    u(x_1,x_2,x_3,t) = t^4 + \sin(x_1)\cdot\sin(x_2) + \sin(x_3)\cdot\sin(t) 
\end{equation}
\end{example}

Within the training set, we have strategically selected 1000 data points for each spacetime under initial conditions, 2000 points spread across the domain, and 500 points allocated to boundary conditions. This deliberate sampling approach ensures that the models undergo extensive training, thereby equipping them with the ability to effectively capture intricate spatiotemporal dynamics. Furthermore, the four models are structured as an assembly of 15 interconnected subnetworks, each corresponding to a specific frequency within the manually defined set $\bm{\Lambda}=(1,2,3,4,5,6,7,8,9,10,11,12,13,14,15)$, thereby enhancing their capacity to learn and represent complex spatiotemporal features with high fidelity. The testing set comprises 13,987 distributed points located within the cut planes, carefully selected while excluding the holes, and encompasses planes that are oriented parallel to both the $xoy$ and $yoz$ coordinate planes, thereby ensuring a thorough evaluation of the model’s performance across multiple spatial dimensions.

\begin{figure}[H]
    \centering
    \subfigure[Analytical solution]{
        \label{E5-01} 
        \includegraphics[scale=0.325]{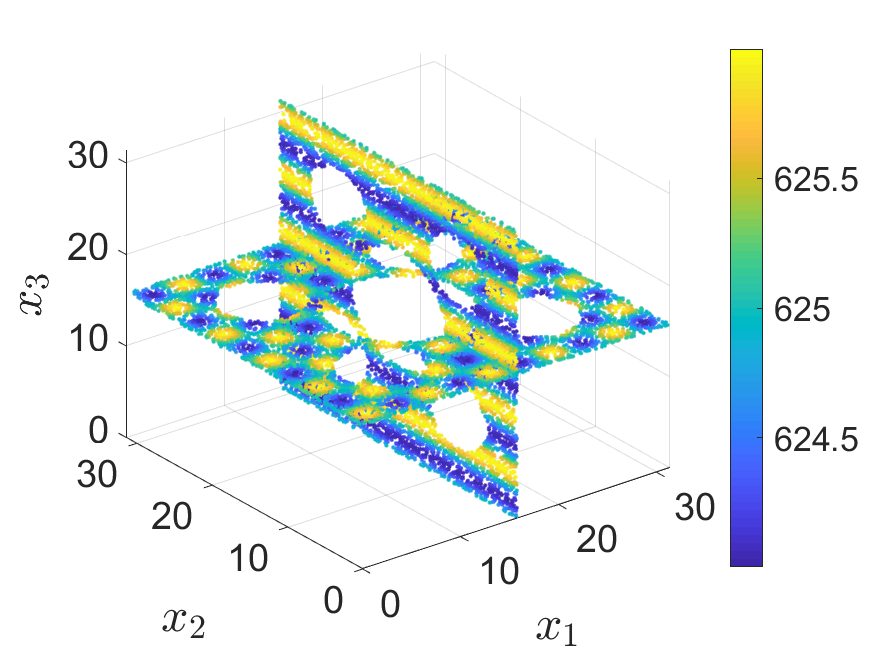}
    }
    \subfigure[Normalized space domain]{
        \label{E5-02} 
        \includegraphics[scale=0.325]{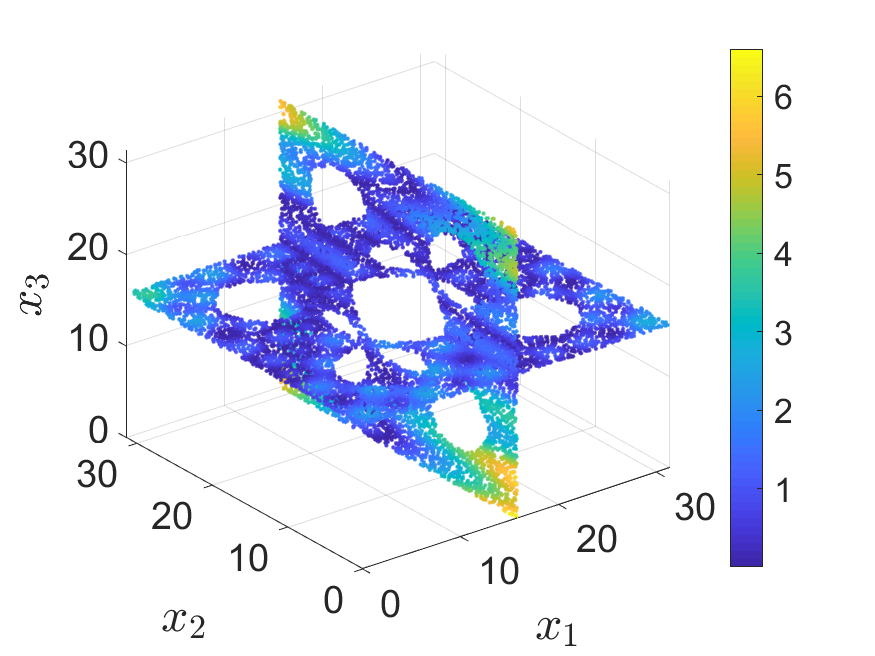}
    }
    \subfigure[Normalized time domain]{
        \label{E5-03}
        \includegraphics[scale=0.325]{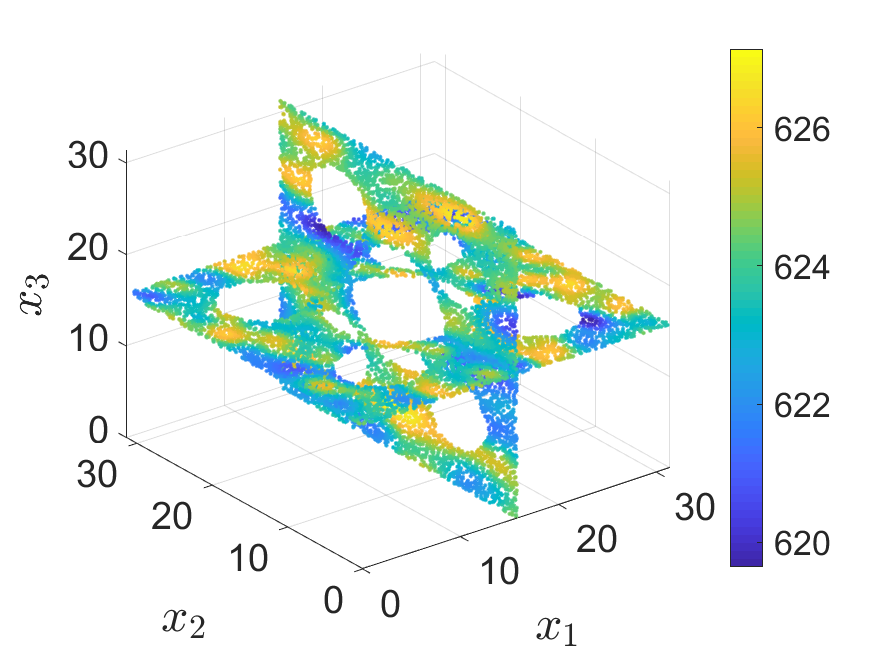}
    }
    \subfigure[Normalized space and time domain]{
        \label{E5-04}
        \includegraphics[scale=0.325]{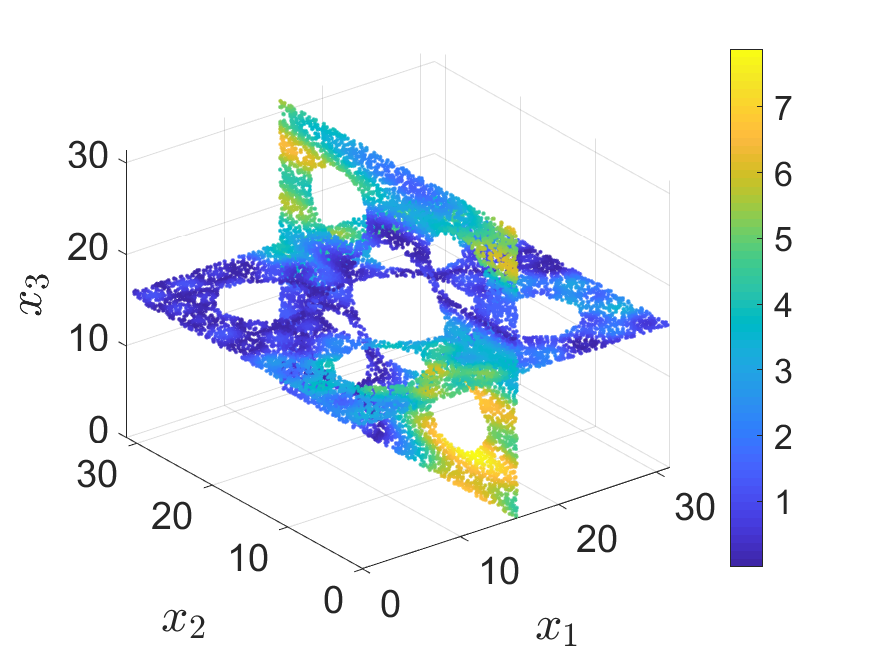}
    }  
    \subfigure[Relative error]{
        \label{E5-05}
        \includegraphics[scale=0.305]{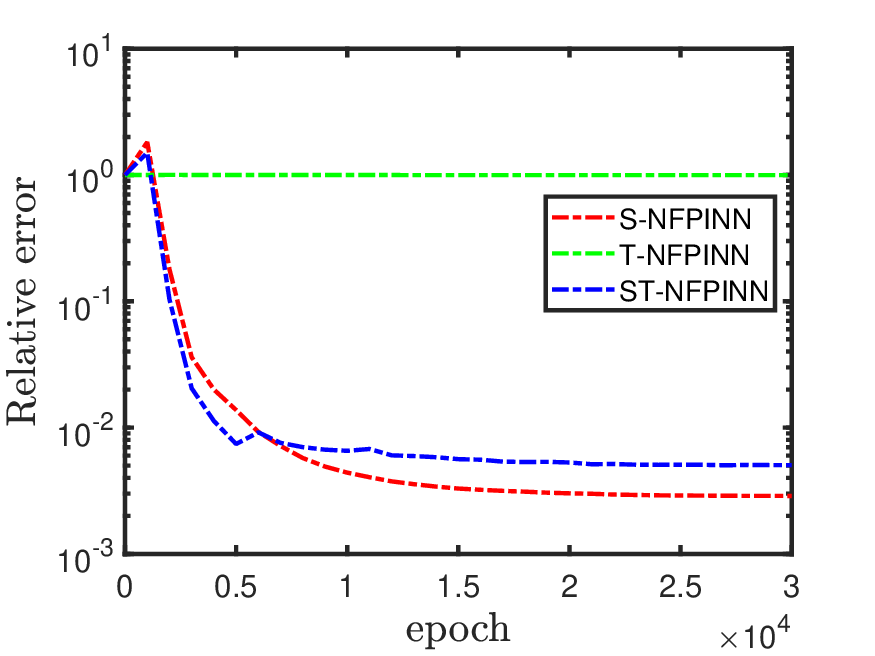}
    }    
    \caption{Numerical results of NFPINN method for test example 5.}
    \label{Fig:Example_5_NFPINN}
\end{figure}

According to the point-wise absolute errors presented in Fig.\ref{E5-01}–Fig.\ref{E5-04}, our proposed S-NFPINN method outperforms T-NFPINN and ST-NFPINN methods in solving Eq.\eqref{eq:Example_5_govern}. The lower point-wise absolute errors distribution demonstrates that the S-NFPINN method effectively captures the underlying wave dynamics with higher accuracy while reducing local discrepancies across the computational domain. Furthermore, the relative errors presented in Fig.~\ref{E5-05} indicate that the S-NFPINN method maintains stability throughout the simulation. These findings confirm that S-NFPINN not only improves solution accuracy but also enhances numerical stability.

\section{Conclusion}\label{sec:6}
To tackle the issue that the classical PINN method and its Fourier enhanced version fail to address wave equations in non-unitized domain and long time range, a normalized Fouirer induced PINN method (NFPINN) is developed to effectively improve the capacity of PINN to deal with wave equations. This new method skilly integrate advantages of Fourier induced PINN and the normalization of space or time cariables. NFPINN is an extension of the FPINN model, which allows the network capturing more spatial and temporal features of wave equations, then the solver of PINN will easily to handle these space-time inputs and product a satisfactory results for wave equations. Furthermore, studies on normalization approaches for spatial and temporal variables have been conducted to determine the optimal normalization method for our proposed approach. Finally, with the aim of verifying the efficiency and robustness for our method, we provide five numerical examples involving both regular and irregular shapes in two-dimensional and three-dimensional Euclidean spaces. The results confirm the accuracy and stability of our proposed method.

\section*{CRediT authorship contribution statement}
Jichao Ma: Methodology, Conceptualization, Validation, Investigation, Writing and Software - Original Draft;
Dandan Liu: Conceptualization, Validation, Investigation and Writing - Original Draft;
Jinran Wu: Writing - Review \& Editing and Project administration;
Xi'an Li: Data curation, Conceptualization, Methodology, Validation, Writing - review $\&$ editing, Funding acquisition.

\section*{Declaration of Competing Interest}
The authors declare that they have no known competing financial interests or personal relationships that could have appeared
to influence the work reported in this paper.

\section*{Acknowledgements}
This paper is funded by the Qingdao Binhai University Research Initiation Fund (Grant No.BS2024A009) and the Natural Science Foundation of Shandong Province, China (Grant No.ZR2024QF057). %The authors are grateful to the anonymous referees for their constructive suggestions and helpful comments, which greatly improved the original manuscript of this paper.

\bibliographystyle{model1-num-names}
\bibliography{References}

\end{document}